\documentclass[12pt]{amsart}
\usepackage{amsfonts}
\usepackage{latexsym}

\newcommand {\supplus}{\mathop{{\supset}\llap{\raise 
0.5pt\hbox{\normalfont\small+}\hskip 0.5pt}}} 

\newcommand {\subplus}{\mathop{{\subset}\llap{\raise 
0.5pt\hbox{\normalfont\small+}\hskip 0.5pt}}}  

\newcommand {\Cee}    {{\mathbb  C}}

\newcommand {\Nee}    {{\mathbb  N}}

\newcommand {\Zee}    {{\mathbb  Z}}

\newcommand {\cal} {\mathcal}

\newcommand {\cX}     {{\cal X}}

\def \opname#1#2%
  {\expandafter\newcommand \csname #1\endcsname {{\mathop{#2}\nolimits}}}

\newcommand{\rmname}[1]
  {\expandafter\newcommand \csname #1\endcsname {{\operatorname{#1}}}}

\newcommand{\rmnameii}[2]
  {\expandafter\newcommand \csname #1\endcsname {{\operatorname{#2}}}}

\rmname{act}
\rmname{Ad}
\rmname{Add}
\rmname{ad}
\rmname{Aff}
\rmname{Alt}
\rmname{alt}
\rmname{Ann}
\rmname{antidiag}
\rmname{Ber}
\rmname{ber}
\rmname{Br}
\rmname{card}
\rmname{ch}
\rmname{Char}
\rmname{cem}
\rmname{cj}
\rmname{Cliff}
\rmname{cntr}
\rmname{codim}
\rmname{coind}
\rmname{const}
\rmname{col}
\rmname{cork}
\rmname{cpr}
\rmname{diag}
\rmnameii{Div}{div}
\rmname{Def}
\rmname{Der}
\rmname{Dim}
\rmname{End}
\rmname{Even}
\rmname{Ev}
\rmname{Ext}
\rmname{gr}
\rmname{Hom}
\rmname{HT}
\rmnameii{Ht}{ht}
\rmname{hwt}
\rmname{Id}
\rmname{id}
\rmname{ind}
\rmname{Ind}
\rmname{Inf}
\rmname{irr}
\rmname{Le}
\rmname{Lie}
\rmname{lwt}
\rmname{mult}
\rmname{Mor}
\rmname{nm}
\rmname{Ob}
\rmname{Odd}
\rmname{ord}
\rmname{Osc}
\rmname{per}
\rmname{Pic}
\rmname{pr}
\rmname{pro}
\rmname{Prime}
\rmname{Proj}
\rmname{prt}
\rmname{pt}
\rmname{Q}
\rmname{qet}
\rmname{qtr}
\rmname{rd}
\rmname{rk}
\rmname{row}
\rmname{Res}
\rmname{salt}
\rmname{Sch}
\rmname{SBr}
\rmname{scalar}
\rmname{sch}
\rmname{sh}
\rmname{Ser}
\rmname{sign}
\rmname{Smbl}
\rmname{spin}
\rmname{ssym}
\rmname{str}
\rmname{st}
\rmname{sgn}
\rmname{sq}
\rmname{symm}
\rmname{supp}
\rmname{Supp}
\rmname{St}
\rmname{Spec}
\rmname{Specm}
\rmname{Spm}
\rmname{tr}
\rmname{vpt}
\rmname{weyl}
\rmname{Weyl}
\rmname{Witt}

\opname{vvol}  {{v\hspace{-0.1ex}o\hspace{-0.02ex}l\/}}
\opname{pnt}  {\text{\normalfont pt}}
\opname{Span} {{Span}}
\opname{slim} {\overline{\lim}}
\opname{Vol}  {{V\hspace{-0.55ex}o\hspace{-0.02ex}l\/}}
\opname{QVol} {{Q\hspace{-0.3ex}V\hspace{-0.55ex}o\hspace{-0.02ex}l\/}}
\opname{PoVol}{{P\hspace{-0.35ex}o\hspace{-0.25ex}V\hspace{-0.55ex}o\hspace{-0.
02ex}l\/}}
\opname{BVol} {{B\hspace{-0.2ex}V\hspace{-0.55ex}o\hspace{-0.02ex}l\/}}
\opname{Par}  {{P\hspace{-0.3ex}a\hspace{-0.05ex}r\/}}

\rmname{Mat}
\rmname{Bil}
\rmname{Diff}
\rmname{Ker}
\rmname{Herm}
\rmname{Coker}
\rmname{Conn}
\rmname{Covect}
\rmname{Vect}
\rmname{Int}

\rmnameii {IM} {Im}
\rmnameii {RE} {Re}

\opname{Aut} {{A\hspace{-0.2ex}u\hspace{-0.1ex}t\/}}
\opname{GL} {{G\hspace{-0.3ex}L}}
\opname{SL} {{S\hspace{-0.3ex}L}}
\opname{Exp} {{E\hspace{-0.2ex}x\hspace{-0.1ex}p\/}}
\opname{G} {{G}}
\opname{GQ} {{G\hspace{-0.2ex}Q}}
\opname{OSp} {{O\hspace{-0.25ex}S\hspace{-0.15ex}p\/}}
\opname{Out} {{O\hspace{-0.25ex}u\hspace{-0.15ex}t\/}}
\opname{Spp} {{S\hspace{-0.2ex}p\/}}
\opname{SpO} {{S\hspace{-0.2ex}p\hspace{-0.02ex}O\/}}
\opname{Pe} {{P\hspace{-0.25ex}e\/}}
\opname{SPe} {{S\hspace{-0.25ex}P\hspace{-0.25ex}e\/}}
\opname{Spin} {{S\hspace{-0.25ex}p\hspace{-0.05ex}i\hspace{-0.1ex}n\/}}
\opname{Iso} {{I\hspace{-0.25ex}s\hspace{-0.1ex}o\/}}
\opname{SSPe} {{S\hspace{-0.25ex}S\hspace{-0.15ex}P\hspace{-0.25ex}e\/}}
\opname{U} {{U\/}}

\opname{cGQ} {{\cal G \hspace{-0.2em} Q \/}}
\opname{cSL} {{\cal S \hspace{-0.2em} L \/}}
\opname{cGL} {{\cal G \hspace{-0.2em} L \/}}
\opname{cOSp} {{\cal O \hspace{-0.2em} S \hspace{-0.3em} \it p\/}}
\opname{cPe} {{\cal P \hspace{-1.5pt} \it e\/}}
\opname{cVect} {{\cal V \hspace{-1.5pt} \it e\hspace{-0.1ex}c\hspace{-0.1ex}t\/
}}
\opname{cVol} {{\cal V \hspace{-1.5pt} \it o\hspace{-0.1ex}l\/}}
\opname{cAut} {{\cal A \hspace{-0.2em} \it u\hspace{-0.1em}t\/}}
\opname{cCovect} {{\cal C \hspace{-1.5pt}
     \it o\hspace{-0.1ex}v\hspace{-0.1ex}e\hspace{-0.1ex}c\hspace{-0.1ex}t\/}}
\opname{CW} {{C\hspace{-0.15ex}W}}

\opname {Ab}   {{\sf Ab}}
\opname {Alg}   {{\sf Alg}}
\opname {ASch}  {{\sf Aff\;Sch}}
\opname {Funct}   {{\sf Funct}}
\opname {Gr}   {{\sf Gr}}
\opname {Grf}  {{{\sf Gr}_f}}
\opname {Mods}   {{\sf Mods}}
\opname {SMods}   {{\sf SMods}}
\opname {Rings}   {{\sf Rings}}
\opname {Salg}   {{\sf Salg}}
\opname {Sets} {{\sf Sets}}
\opname {SSMan} {{\sf SMan}}
\opname {Top}  {{\sf Top}}
\opname {Vebun}   {{\sf Vebun}}

\newcommand {\tto} {\longrightarrow}
\newcommand {\pder}[1] {{\frac{\partial}{\partial {#1}}}}
\newcommand {\pderf}[2] {{\frac{\partial {#1}}{\partial {#2}}}}

\newcommand {\bcdot}   {\mathbin{\hbox{\raise.4ex\hbox{\bf.}}}} 

\newcommand {\secno} {}
\newcommand {\ssecfont} {\normalfont\bf}

\newtheorem{Theorem}{\secno Theorem}

\newtheorem{Lemma}[Theorem]{\secno Lemma}

\newtheorem{Corollary}[Theorem]{\secno Corollary}

\newenvironment {th*}[1]
    {\gdef\thname{#1} \begin{thn}}%
    {\end{thn}}
\newtheorem{thn}[Theorem] {\thname}

\theoremstyle{definition}

\newenvironment {ex*}[1]
    {\gdef\thname{#1} \begin{exn}}%
    {\end{exn}}
\newtheorem{exn}[Theorem]{\thname}

\theoremstyle{remark}

\newtheorem{Remark}[Theorem]{\secno Remark}

\newenvironment {rem*}[1]
    {\gdef\thname{#1} \begin{remn}}%
    {\end{remn}}
\newtheorem{remn}[Theorem]{\thname}

\newcommand {\ssec}{\subsection*}
\newcommand {\ssbegin}[2]
  {\def \secno {\gdef \secno {}{\ssecfont #1. }}%
   \begin{#2}}
\begin{document}

\title{Invariant functions on supermatrices}

\author{Vladimir Shander}

\address{c/o D. Leites, Dept. of Math., Univ. of Stockholm, Roslagsv.
101, Kr\"aftriket hus 6, S-106 91, Stockholm, Sweden; 
mleites@matematik.su.se}


\begin{abstract} There are two superanalogs of the general linear 
group: $\GL (m|n)$ and $\GQ(n)$; for any supercommutative superalgebra 
$\Lambda$ let $\G(\Lambda)$ be the set of $\Lambda$-points of the 
supermanifold $\G$.  Here there are described the $\GQ (n; 
\Lambda)$-invariant functions on $\Q(n; \Lambda )$ and $\GL(n; \Lambda 
)$-invariant functions on $\Mat (n|n; \Lambda )_{\bar 1}$, more 
exactly, the invariants of the action of Lie supergroups $\GQ(n)$ and 
$\GL (n|n)$ on the supermanifolds corresponding to $\Q(n; Lambda )$ 
and $\Mat (n|n; \Lambda )_{\bar 1}$.

The obtained answer is interpreted in terms of $\GL (n; \Lambda 
)$-invariants on $\Q(n; \Lambda )$ and $\GL (n|n; \Lambda 
)$-invariants on $\Mat (n|n; \Lambda )_{\bar 1}$ described 
functorially in (i.e., independently of) $\Lambda$.
\end{abstract} 

\thanks{I am thankful to D.~Leites, A.~Levin, A.~Premet, V.~Retakh 
A.~Sergeev, and A.~Vaintrob for useful discussions and to S.~Lvovsky 
who pointed out at the paper \cite{Se} and made elucidating remarks.  To 
D.~Leites and the Swedish NFR I am thankful for hospitality and 
financial support during my stay in Sweden in the autumn of 1990.}

\section* {Introduction} 

This is an expanded version of my paper {\it Orbits and invariants of 
the supergroup ${\rm GQ}\sb n$}.  Funktsional.  Anal.  i Prilozhen.  
26 (1992), no.  1, 69--71 (in Russian)

\ssec{Problem formulation} For classical Lie 
groups the problem of description of invariants is completely solved 
by H.~Weyl.  In the theory of (super)matrices over a supercommutative 
superalgebra $\Lambda $ there naturally arise several problems of 
description of invariants (the necessary background is given in \S 1):
 
There are two superanalogs of the matrix algebra.  They are denoted by 
$\Mat (p|q; \Lambda )$ and $\Q(n; \Lambda )$, respectively, on each of 
these analogs the corresponding group of invertible matrices, $\GL 
(p|q; \Lambda )$ or $\GQ(n; \Lambda )$, respectively, acts by 
conjugations.  The actions of $\GL (p|q; \Lambda )$ and $\GQ(n; 
\Lambda )$ preserve the parity of matrices, i.e., the decomposition of 
matrices into even and odd ones:
$$
\Mat (p|q; \Lambda )_{\bar 0}\oplus \Mat (p|q; \Lambda )_{\bar 1}\; \; 
\text{ and }\Q(n; \Lambda )=\Q(n; \Lambda )_{\bar 0}\oplus \Q(n; 
\Lambda )_{\bar 1}
$$
and, therefore, the invariants of even matrices and the invariants of 
odd matrices should be described separately.  Taking into account a 
canonical isomorphism $\Q(n; \Lambda )_{\bar 0}\cong \Q(n; \Lambda 
)_{\bar 1}$ of $\GQ (n; \Lambda )$-modules, we get the following three 
problems of description of invariants of supermatrices:

1) $\GL (p|q; \Lambda )$-invariant functions on $\Mat (p|q; \Lambda 
)_{\bar 0}$;

2) $\GL (p|q; \Lambda )$-invariant functions on $\Mat 
(p|q; \Lambda )_{\bar 1}$; 

3) $\GQ(n; \Lambda )$-invariant functions on $\Q(n; \Lambda )_{0}$.

These are precisely the problems that we are going to study.  More 
exactly, in accordance with the general principles stated in \cite{L}, v.  
30, \S 2.4, we should consider these problems functorially in 
$\Lambda$, which in simpler words means that the answer should not 
depend on $\Lambda $.

Recall (cf. \cite{W}) that the description of invariants of $\GL (n; \Cee )$
consists of two separate statements:

$\bullet$ first, any invariant function of a matrix $A\in \Mat (n; \Cee )$
is a function in $\tr A, \dots, \tr A^{n}$ and

$\bullet$ second, polynomial invariant functions {\it polynomially} 
depend on $\tr A, \dots, \tr A^{n}$.

The first statement is related to the fact that in $\Mat (n; \Cee )$ 
there exists a dense set of matrices that can be reduced to the 
diagonal form and the second statement is related with the theorem on 
symmetric polynomials.

In the classical invariant theory the algebra of invariant polynomials 
usually has a finite set of polynomial generators (N\"otherian 
property).  In the supercase, as well as over fields of prime 
characteristic, this is not so and the above listed problems give us 
counterexamples.

For the first of these problems, however, F.~Berezin and Kac  
proved that any invariant polynomial on $\Mat (p|q; \Cee )_{\bar 0}$ 
can be expressed in terms of $p+q$ polynomials in $\str A, \dots , 
\str A^{p+q}$ {\it but not necessarily in a polynomial way}, for 
instance, as a ratio of two polynomials.

This theorem whose proof is only contained in the English version of 
Berezin's posthumous book \cite{Be} (see also \cite{Ka}) also 
consists of two statements:

(1) on existence of a dense in $\Mat (p|q; \Cee )_{\bar 0}$ set of 
diagonalizable matrices with pair-wise distinct eigenvalues and

(2) on a possibility to {\it rationally} express any polynomial in 
(even) variables $\lambda _{1}, \dots , \lambda _{p}, \mu _{1}, \dots 
, \mu _{q}$ symmetric, separately, in $\lambda $ and in $\mu $ in 
terms of functions $s_{1}(\lambda , \mu ), \dots , s_{p+q}(\lambda , 
\mu )$, where
$$
s_{i}(\lambda , \mu )=\sum_j \lambda ^{i}_{j}-\sum_j \mu ^{i}_{j}.
$$

Taking into account the fact that the remarkable invariant function on 
$\Mat (p|q; \Cee )_{\bar 0}$ --- the Berezinian --- is not polynomial 
but is a rational function, it is natural to interpret the above 
mentioned theorem as the following statement {\sl the algebra of 
invariant rational functions on $\Mat (p|q; \Cee )_{\bar 0}$ is 
isomorphic to the algebra of rational functions in $p+q$ even 
generators $\str A, \dots , \str A^{p+q}$}, cf.  \cite{Ka}.

\ssec{Related results} At present, the knowledge of 
$\GQ(n)$-invariants is scanty.  On $\Q(n)_{\bar 0}$, there are 
analogues of the trace and determinant --- odd $\GQ(n)$-invariant 
functions, $\qtr A$ and $\qet A$, the first of which is linear and the 
second one rational, \cite{BL}.

A.~Sergeev proved in \cite{S} that the algebra of
$\GQ$-invariant {\it polynomials} on $\Q(n)_{\bar 0}$ is generated by the
infinite set of polynomials
$$
\qtr  A, \qtr  A^{2}, \dots
, \qtr  A^{n}, \dots  .
$$
An implicit description of invariant polynomials on $\Q(n)$ is also 
contained in \cite{F}. ??

There are natural maps
$$
\Mat (p|q; \Lambda )_{\bar 1}\tto \Mat (p|q; 
\Lambda )_{\bar 0}; \; \; \; \; M\mapsto M^{2} 
$$
and
$$
\Q(n, \Lambda )_{\bar 0}\tto \Mat (n|n, 
\Lambda )_{\bar 0};\; \; \; \; (A_{0}+A_{1})\mapsto 
\begin{pmatrix} A_{0} & A_{1} \cr A_{1} & A_{0}\end{pmatrix} 
$$
but these maps do not give any invariants for $\Mat (p|q; \Lambda 
)_{\bar 1}$ and $\Q(n, \Lambda )$, since the supertraces of the images 
under the above mappings are zero.

Observe two discouraging circumstances concerning invariant
functions on $\Q(n)$.

First, one can show (see \S 3) that there is no
finite set of invariant functions that can generate all
the other invariant functions. Second, the invariant functions on $Q(n)$
carry very few information on the corresponding $\GQ
(n)$-orbits. In particular, it is impossible to determine
from the value of all the invariant polynomials at a
given matrix whether the matrix is invertible or not.

\ssec{Our result} Here I describe $\GQ (n; \Lambda)$-invariant 
functions on $\Q(n; \Lambda )$ and $\GL(n; \Lambda )$-invariant 
functions on $\Mat (n|n; \Lambda )_{\bar 1}$ (functorially in 
$\Lambda$).  More exactly, I describe {\it the invariants of the 
action of Lie supergroups $\GQ(n)$ and $\GL (n|n)$ on the 
corresponding supermanifolds denoted by $\Q(n)$ and $\Odd (n)$.}

The answer obtained is interpreted in terms of $\GL (n; \Lambda 
)$-invariants on $\Q(n; \Lambda )$ and $\GL (n|n; \Lambda 
)$-invariants on $\Mat (n|n; \Lambda )_{\bar 1}$ described 
functorially in $\Lambda$.

In what follows, in order to describe the invariant functions we will 
use {\it semi-invariants}, the functions which are not invariant but 
which under the action of the supergroup accrew summands that belong 
to the ideal generated by invariant functions.

It turns out that any $\GQ(n)$-invariant function on $Q(n)$ can be 
expressed as a function in $n$ odd invariants $\qtr A, \dots , \qtr 
A^{n}$ and $n$ {\it odd noninvariant rational functions}
$$
t_{1}(A), \dots , t_{n}(A).
$$

We will start with holomorphic invariant functions and then
pass to the rational and polynomial functions.

The algebra of $\GL (n|n)$-invariant functions on $\Odd (n)$ turns out 
to be isomorphic to the algebra of $\GQ (n)$-invariant functions (of 
the same class, i.e., holonomic, rational or polynomial functions, 
respectively) on $\Q(n)$ in spite of the fact that this isomorphism of 
algebras of invariants is not induced by any natural map of $\Q(n)$ to 
$\Odd (n)$, or the other way round.

The description of all invariant functions on odd matrices of general 
form is also similar to that of invariants on $\Q(n)$ but is slightly 
more cumbersome and will be given separately.

The contents of the paper is as follows.

\S 1 contains the necessary background.  In \S 2 we consider in detail 
the case $n=1$.  \S 3 is devoted to the study of functions on $\Cee 
^{n|n}$ invariant with respect to the action of the symmetric group 
$S_{n}$ and $0|n$-dimensional abelian supergroup $\Cee ^{0|n}$.  In 
particular, I prove a superanalog of the theorem on symmetric 
functions:

{\sl Any symmetric function in $n$ non-homogeneous variants can be 
uniquely expressed in terms of $n$ symmetric non-homogeneous {\em 
(with respect to parity)} functions}.

In \S 4 I prove the main results on invariant
functions on $\Q(n)$ and $\Odd (n)$.  \S 5 contains several 
examples.

\section* {\protect \S1. Background}

\ssec {1.0} The sources of information on superalgebra and 
supercalculus are \cite{Be}, \cite{L}, \cite{Ma}.  Recall some 
notations and definitions.  If $C=C_{\bar 0}\oplus C_{\bar 1}$ is a 
supercommutative superalgebra and $I_{C}$ the ideal in $C$ generated 
by $C_{\bar 1}$, then the image of $x\in C$ in $C/I_{C}$ is denoted by 
$\cpr x$.

At first, let $C$ be an arbitrary algebra.  Then the algebra of 
$n\times n$-matrices $\Mat (n; C)$ acts from the left on the columns 
of length $n$.  If $L$ is a free right $C$-module of $\rk~ n$, then, 
having fixed a basis $e_{1}, \dots , e_{n}$ in $C$, we determine an 
isomorphism $\Hom _{C}(L, L)\simeq \Mat (n; C)$ such that the matrix 
of each operator acts from the left on the column of right 
coordinates; the change of basis is performed with the help of 
invertible matrix --- an element from $\GL (n; C)$ --- and acts in the 
usual way on the matrices of operators.  If an element $c\in C$ 
belongs to the center of $C$, then its action on $L$ is given by a 
scalar matrix.

\ssec {1.1} If $C$ is a supercommutative superalgebra, then $\Mat (n; 
C)$ and $\GL (n; C)$ are also denoted by $\Q(n; C)_{\bar 0}$ and 
$\GQ(n; C)$, respectively; for motivations see \cite{L}, v.  30, Ch.  
1.  Any matrix $A\in \Q(n; C)_{\bar 0}$ can be uniquely expressed in 
the form $A_{0}+A_{1}$, where all the elements of $A_{0}$ are even, 
all the elements of $A_{1}$ are odd, and $\cpr A=\cpr A_{0}\in \Mat 
(n; C/I_{C})$.

A matrix $A$ is invertible if and only if the matrix $\cpr A$ is 
invertible and $A$ is nilpotent if and only if $\cpr A$ is nilpotent.

In particular, $\GQ(n; C)$ consists of all the matrices
$A_{0}+A_{1}$ such that $\cpr  A_{0}$ is invertible.

The analogues of the trace and determinant for $\Q(n; C)$ are
the $\GQ(n)$-invariant functions $\qtr  : \Q(n; C)_{\bar 0}\tto
C_{\bar 1}$ and $\qet : \GQ(n)\tto C_{\bar 1}$, where
\begin{align*}
\qtr  (A_{0}+A_{1})&=\tr  A_{1}, \cr
\qet (A_{0}+A_{1})&=\sum^{}_{1\le i}(\frac{1}{i})\tr  
(A^{-1}_{0}A_{1})^{i}. \cr
\end{align*}

\ssec{1.2} On a free right $C$-module $L$ with a basis
$e_{1}, \dots
, e_{p+q}$ introduce a parity $(\Zee _{2}$-grading) declaring the
first $p$ elements of the basis even and the other ones odd;
let the $C$-action on $L$
be an even map $L\times C\tto L$. Then
$$
\Hom  (L, L)=\Hom_{\bar 0}(L, L)\oplus \Hom_{\bar 1}(L, 
L), 
$$
 where the operators from $\Hom_{\bar 0}(L, L)$ 
preserve and the
operators from $\Hom_{\bar 1}(L, L)$ change the parity of
homogeneous elements of $L$.

The parity of the elements from the space $\Mat (p+q; C)$ is
defined as follows: the matrices are split into blocks
$$
\cX= \begin{pmatrix} X & Y \cr Z & T \cr\end{pmatrix}, \; \text{ where $X\in 
\Mat
(p; C)$, $T\in \Mat (q; C)$}, \eqno{(*)}
$$
and the matrix $\cX$ is even if the elements
of the blocks $X$ and $T$ are even and the elements of the blocks $Y$
and $Z$ are odd whereas $\cX$ is odd if the
elements of $X$ and $T$ are odd and the elements of $Y$ and $Z$ are even.

The algebra $\Mat (p+q; C)$ with the above parity is denoted by $\Mat 
(p|q; C)$ and the group of even invertible matrices is denoted by $\GL 
(p|q; C)$; it is identified with the changes of basis that preserve 
the parity and the {\it format} of the matrix (see \cite{Ma} or 
\cite{L}, v.  30, Ch.1).

The {\it supertrace} $\str:\Mat (p+q; C)\tto C$ is 
determined on the even matrices $(*)$ as $\tr X-\tr T$ and on odd 
matrices as $\tr X+\tr T$.  The supertrace is $\GL (p|q)$-invariant.

\ssec {1.3} The following Theorem was proved in 1978 but
was not published. Since it is not covered by subsequent
publications --- only in \cite{Be} the first two headings of
Corollary are proved --- we give a complete proof, 
especially since it is so simple.

\begin{Theorem} Let $C$ be a local supercommutative
superalgebra with unit over an algebraically closed field
$k$ and all the elements of the maximal ideal $I\subset C$ are
nilpotent. Then:

$1)$ For any free $C$-module $L$ of $\rk~ \; n$ and any
$A\in \Hom  _{C}(L, L)$ the submodules $L(\lambda )=
\bigcup\limits _{i}{\rm Ker} (A-\lambda )^{i}, \lambda \in C$, are
either zero or free and $L=\bigoplus\limits _{\lambda \in k}L(\lambda )$.

$2)$ If $A\in \Q(n; C)_{\bar 0}$ and the collection of distinct eigenvalues of
$\cpr  A\in \Mat (n; k)$ {\em (we apply the canonical
projection $\cpr : C\tto C/I$ to a matrix element-wise)} is equal to
$\lambda _{1}, \dots , \lambda _{s}$ with multiplicities $n_{1}, \dots , 
n_{s}$,
respectively, then there exists a matrix
$G\in \GQ (n; C)$ such that
$G^{-1}AG$ is of a block-diagonal form $\begin{pmatrix} A_{1} & \dots & 0 \cr 
\vdots& \vdots& \vdots&\cr 0 & \dots & A_{s} \cr\end{pmatrix}$, where 
$A_{i}\in Q
(n_{i}; C)$ and $\lambda _{i}$ is the only eigenvalue of the  matrix $\cpr 
A_{i}$. The matrix $A_{i}$ corresponding to the eigenvalue $\lambda _{i}$ is
uniquely defined up to the $\GL (n_{i}; C)$-action.

$3)$ If $A=\begin{pmatrix} X & Y \cr Z & T \cr\end{pmatrix}\in \Mat  (p|q;
C)_{\bar 0}$ and the collection of distinct eigenvalues of
$\cpr  A=\begin{pmatrix} \cpr  X & 0 \cr 0 & \cpr   T \cr\end{pmatrix}$ is 
equal
to $\lambda _{1}, \dots , \lambda _{s}$, then there exist elements
$G\in \GL (p|q; C)$ and a partition of the set $\{1, \dots , p+q\}$ into
nonintersecting subsetes $J_{1}, \dots , J_{s}$ such that all the elements of
$G^{-1}AG$ for which the number of the row and the number of the column belong
to the distinct subsets of the partition are zeros. 

The square submatrix $A_{i}$
corresponding to $J_{i}$ belongs to $\Mat (p_{i}| q_{i}; C)_{\bar 0}$, where
$p_{i}$ and $q_{i}$ are the multiplicities of $\lambda _{i}$ in the spectra of
the matrices $\cpr  X$ and $\cpr  T$, respectively, and  where $\lambda _{i}$ 
is
the only eigenvalue of $\cpr  A_{i}$. The matrix $A_{i}$  corresponding to
$\lambda _{i}$ is determined uniquely up to the action of $\GL (p_{i}|q_{i};
C)$.

$4)$ If $A=\begin{pmatrix} X & Y \cr Z & T \cr\end{pmatrix}\in \Mat (p|q; 
C)_{\bar 1}$ and the set of eigenvalues of
$\cpr  A^{2}$ is equal to $\lambda _{1}, \dots , \lambda _{s}$, then there 
exist
a matrix $G\in \GL (p|q; C)$ and a partition of the set $\{1, \dots , p+q\}$ 
into
nonintersecting subsets $J_{1}, \dots , J_{s}$ such that all the elements of
$G^{-1}AG$ for which the number of the row and the number of the column belong
to distinct subsets of the partition vanish and the square matrix $A_{i}$
corresponding to $J_{i}$ belongs to $\Mat (p_{i}|q_{i}; C)_{\bar 1}$,  where
$p_{i}$ and $q_{i}$ are the multiplicities of $\lambda _{i}$ in the spectra of
$\cpr (YZ)$  and $\cpr  (ZY)$, respectively, and $\lambda _{i}$ is the only
eigenvalue of $\cpr   A^{2}_{i}$. 

The matrix $A_{i}$ corresponding to $\lambda _{i}$ is determined uniquely up to
the action of $\GL (p_{i}|q_{i}; C)$ and if $\lambda _{i}\neq 0$, then
$p_{i}=q_{i}$ and $A_{i}$ can be reduced to the form 
$\begin{pmatrix} R & T \cr 1_n & 0 \cr\end{pmatrix}$,  where $\cpr  T$ has an
only eigenvalue,  $\lambda _{i}$.
\end{Theorem}

\begin{Corollary} Under the same conditions on $C$ as in Theorem:

$1)$ If $A\in \Q(n; C)_{\bar 0}$ and $\cpr  A$ have no multiple
eigenvalues, then the $\GQ(n; C)$-action can reduce $A$ to a diagonal
form.

$2)$ If $A\in \Mat (p|q; C)_{\bar 0}$ and $\cpr  A$  have no multiple
eigenvalues, then the $\GL (p|q; C)$-action can reduce $A$ to a
diagonal form.

$3)$ If $A\in \Mat (n|n; C)_{\bar 1}$ and $\cpr   A^{2}$ has no multiple
eigenvalues, then the $\GL (n|n; C)$-action can reduce $A$ to the
form
$\begin{pmatrix} R & T \cr 1_n & 0 \cr\end{pmatrix}$, where $R$ and 
$T$ are diagonal matrices.
\end{Corollary}

\begin{Remark} 1) In this work the theorem is only used in
the case when $C$ is the Grassmann (exterior) algebra over $\Cee $
and $I=I_{C}$. A more general formulation guarantees
the possibility to work, if needed, with algebras obtained
from the Grassmann algebra by adjoining to it roots of
algebraic equations. For example, if $C$ is a
supercommutative superalgebra over $\Cee $ with two odd
generators $\xi _{1}$ and $\xi _{2}$, one even generator $t$ satisfying the
relation
$t^{2}=\xi _{1}\xi _{2}$, then $C/I_{C}$ is a two-dimensional algebra but,
nevertheless, $C/I=\Cee $ and Theorem is applicable to matrices
over $C$.

2) The formulations for $A\in \Mat (p|q; C)$ are cumbersome
because we have restricted ourselves to matrices of the standard format.
If we allow arbitrary farmats, then we can
reduce matrices from $\Mat (p|q; C)$ to the conventional 
Jordan block-diagonal form.
\end{Remark}

{\bf Proof of Theorem}. The first two statements are
proved simultaneously. Fix a basis in $L$ and denote by
$A^{(0)}$ the matrix of the operator $A$ in this basis. Since
$C/I=k$, then there exists a finite dimensional subalgebra
of $C$ containing all the elements of $A^{(0)}$ and, therefore, in
what follows we may assume that the decreasing filtration
of $C$ with respect to powers of $I$ is finite.

Let us prove that if $\cpr  A^{(0)}=\begin{pmatrix} B & 0 \cr 0 & C
\cr\end{pmatrix}$, where $B\in \Mat
(n_{1}; k)$ and $D\in \Mat (n_{2}; k)$ have no common eigenvalues, then
for any $i>0$ there exists $G_{i}\in \GL (n; C)$ such that $\cpr  
G_{i}=1$
and
$$
A^{(i)}=G^{-1}_{i}A^{(0)}G_{i}\equiv \begin{pmatrix} B_{i} & 0 \cr 0 & D_{i}
\cr\end{pmatrix}(\mod  I^{i}).
$$
Let $A^{(i)}\equiv \begin{pmatrix} B_{i} & 0 \cr 
0 & D_{i} 
\cr\end{pmatrix} (\mod I^{i})$.  Then $A^{(i)}\equiv \begin{pmatrix} 
B_{i}+\Delta _{1} & \Delta _{2} \cr 
\Delta _{3} & D_{i}+\Delta 
\cr\end{pmatrix}$, where $\Delta _{i}\equiv 0 (\mod I^{i})$.  Since 
$B$ and $D$ have no common eigenvalues, 
then the linear maps $x\mapsto 
Bx-xD$ and $y\mapsto Dy-yB$ defined on the spaces of $n_{1}\times 
n_{2}$-matrices and $n_{2}\times n_{1}$-matrices with elements from 
the ground field $k$, respectively, are one-to-one (\cite{G}, Ch.  8).  
Therefore, there exists a matrix $\Delta \in I^{i}\cdot \Mat (n; C)$ 
such that $\left[\begin{pmatrix} B & 0 \cr 
0 & D \cr\end{pmatrix}, 
\Delta \right]=-\begin{pmatrix} 0 & \Delta _{2} \cr 
\Delta _{3} & 0 
\cr\end{pmatrix}$; hence, $1+ \Delta \in \GQ(n; C)$ and
$$
(1+\Delta )^{-1}A^{(i)}(1+\Delta )\equiv A^{(i)}+[A^{(i)}, \Delta ]\equiv 
\begin{pmatrix} B_{i}+\Delta _{1} & 0 \cr 0 & C_{i}+\Delta _{4}
\cr\end{pmatrix} (\mod  I^{i+1}).
$$
Thus, having first reduced $\cpr  A^{(0)}$ to the Jordan form,
we get a basis of $L$ in which the matrix of the operator $A$ is
of block-diagonal form described in heading 2) of
Theorem: $A^{(\infty )}=\begin{pmatrix} A_{1} & \dots & 0 \cr 
\vdots&\vdots&\vdots \cr
0 & \dots & A_{s} \cr\end{pmatrix}$.

Since $\cpr  A_{i}$ only has one eigenvalue, $\lambda _{i}$, then $A_{i}-
\lambda $ is
invertible for $\lambda \neq \lambda _{i}$ and $A_{i}-\lambda _{i}$ is 
nilpotent. This means
that the submodule of $L$ spanned by basis vectors
corresponding to the block $A_{i}$ coincides with $L(\lambda _{i})$ and if
$\lambda $ does not coincide with any of the eigenvalues of $\cpr  
A^{(0)}$, then $L(\lambda )=0$. Now, the first two heading of Theorem
are completely proved: the matrices $A_{i}$ are uniquely
determined up to a choice of a basis in $L(\lambda _{i})$.

3) If $A$ is an even operator, then all the $L(\lambda _{i})$ are
homogeneous submodules and selecting in each of them a
basis consisting of homogeneous elements we get 
heading 3).

4) If $A$ is an odd operator, then the $L(\lambda )$ are not
homogeneous submodules. But $A^{2}$ is an even operator and heading 3) is
applicable to it.

If $A_{i}=\begin{pmatrix} X & Y \cr Z & T \cr\end{pmatrix} \in \Mat (p_{i}|
q_{i}; C)_{\bar 1}$ and $\cpr  A^{2}_{i}=\begin{pmatrix} \cpr  
YZ & 0 \cr 0 &
\cpr  ZY \cr \end{pmatrix}$ has only one eigenvalue $\lambda _{i}\neq 0$, 
then $p_{i}=q_{i}$ and $Y$ and $Z$
are invertible. Hence,
$$
\begin{pmatrix} Z & T \cr 0 & 1 \cr\end{pmatrix}\begin{pmatrix} X & Y \cr Z & T
\cr\end{pmatrix}
\begin{pmatrix} Z & T \cr 0 & 1 \cr\end{pmatrix}^{-1}=
\begin{pmatrix} T+ZXZ^{-1} & ZY-ZXZ^{-1}T \cr 1 & 0 \cr\end{pmatrix}
$$
 is of the form  desired. Theorem is proved. \qed

\ssec {1.4} Let $M$ be a supermanifold. Denote by $F_{M}$
or just by $F$ the sheaf of functions on $M$ and if $f$ is a
function on $M$, i.e., a (global) section of $F$, then $\cpr  f$
is identified with the restriction of $f$ onto a
canonically embedded into $M$ underlying manifold denoted
by $M_{rd}$. An open subsupermanifold $V\subset M$ is determined, see 
\cite{L}, v.
30, Ch. 3,  by the open subset $V_{rd}\subset M_{rd}$. 

In what follows we will work with complex-analytic supermanifolds (cf.  
\cite{Ma}) and, except for subsection 4.6, all supermanifolds are 
superdomains, i.e., open subsuperdomains in $\Cee ^{p|q}$.  This means 
that their underlying manifolds are domains in $\Cee ^{p}$ and $F$ is 
the sheaf of analytic functions on $\Cee ^{p}$ with values in the 
Grassmann algebra with $q$ indeterminates.

On $\Cee ^{p|q}$, there exists a global coordinate system consisting 
of $p$ even functions $u_{1}, \dots , u_{p}$ and $q$ odd functions 
$\xi _{1}, \dots , \xi _{q}$.  An arbitrary function $f\in F(M)$ can 
be uniquely expressed in the form $f=\sum f_{\alpha }(u_{1}, \dots , 
u_{p})\xi ^{\alpha _{1}}_{1}\dots \xi ^{\alpha _{q}}_{q}$, where 
$\alpha $ runs over $\{1, 0\}^{q}$.  The morphism of superdomains 
$\varphi : V\tto W\subseteq C^{p|q}$ is determined by the 
morphism of superalgebras $\varphi ^{*}:F(W)\tto F(V)$ 
which in turn is uniquely defined by its coordinate expression --- the 
collection of $p$ even and $q$ odd functions $\varphi ^{*}(u_{1}), 
\dots , \varphi ^{*}(u_{p}), \varphi ^{*}(\xi _{1}), \dots , \varphi 
^{*}(\xi _{q})$.

We will only need supermanifolds associated with
$$
\Mat (p|q; \Lambda )_{\bar 0}, \; \; \Mat (p|q; 
\Lambda )_{\bar 1}, \; \; \GL (p|q; \Lambda ), \; \; \Q(n; 
\Lambda )_{\bar 0}.
$$
Intentional similarity of notations when we deal with distinct categories will
not cause a misunderstanding since it is always clear from the contents which
category we are talking about.

It is convenient to think that the coordinates on
the supermanifold $\Q(n)=\Cee ^{n^{2}|n^{2}}$ fill in two square
matrices of size $n\times n$ each: $X=(X_{ij})$ that consists of even
coordinates and $\xi =(\xi _{ij})$ that consists of odd coordinates. Clearly, 
$Q
(n)_{rd}=\Mat (n; \Cee )$; the supermanifold $\GQ(n)$ is an open
subsupermanifold in $\Q(n)$ and the underlying group of 
$\GQ(n)$ is $\GL (n; \Cee )$.

The action $\ad : \GQ(n)\times \Q(n)\tto \Q(n)$ is 
defined which
in coordinates $X, \xi $ on $\Q(n)$ and $Y, \eta $ on $\GQ(n)$ is given by
the formula
$$
\ad ^{*}(X+\xi )=(Y+\eta )^{-1}(X+\xi )(Y+\eta ).
$$

Similarly, the coordinates on the supermanifold of even
matrices $\Ev  (p|q)=\Cee ^{p^{2}+q^{2}|2pq}$ fill 
out the even matrix $\begin{pmatrix} X & Y \cr Z &
T \cr \end{pmatrix}$, where the matrices $X$ and $T$ are filled out by even
coordinates whereas the elements of $Y$ and $Z$ are odd
coordinates. The supergroup $\GL (p|q)$ is an open subsupermanifold of
$\Ev  (p|q)$ such that
$$
\GL (p|q)_{rd}=\Ev  (p|q)_{rd}\cap \GL (p+q; \Cee ).
$$
The coordinates on the supermanifold $\Odd
(p|q)=\Cee ^{2pq|p^{2}+q^{2}}$ of odd matrices fill out the 
odd matrix $\begin{pmatrix} X'&
Y' \cr Z' & T' \cr\end{pmatrix}$ and the coordinate expression of the action
$$
\ad  : \GL (p|q)\times \Odd (p|q)\tto 
\Odd (p|q)
$$
is similar to the action $\ad $ of $\Q(n)$:
$$
\ad  ^{*}\begin{pmatrix} X' & Y' \cr Z' & T' \cr\end{pmatrix}=
\begin{pmatrix} X & Y \cr Z & T \cr\end{pmatrix}^{-1}\begin{pmatrix} X'& Y'\cr
Z'& T'\cr\end{pmatrix}
\begin{pmatrix} X & Y \cr Z & T \cr\end{pmatrix}.
$$
In what follows the notation $\Odd (n|n)$ is abbreviated to $\Odd (n)$.

\ssec{1.6} In contradistinction to the classical calculus, in supercalculus
the function on $M$ is not defined by its values at $\Cee $-points
of $M$ and, therefore, one has to explicitly introduce
dependence on parameters, cf. \S 2.5.

For convenience, we assume that parameters run over an arbitrary 
supermanifold $U$ though it suffices to take as $U$ \lq\lq purely odd" 
supermanifolds $\Cee ^{0|s}$ with sufficiently large $s$.  If $M$ and 
$U$ are supermanifolds, then a $U$-{\it family of points of} $M$ is 
any morphism $\varphi :U\tto M$ and the function $f$ on 
$U\times M$ is called a $U$-{\it family of functions on} $M$, etc.  
The necessity to introduce parameters and the corresponding technique 
is discussed in detail in \cite{L}.  When the work with parameters can 
be performed automatically we will not mention them.

\section* {\protect \S 2. Invariant functions on $\Q(1)$ and $\Odd(1)$}

\ssec{2.1} The action $\rho : \GQ(1)\times \Q(1)
\tto \Q(1)$ of the
supergroup $\GQ(1)$ is given in the standard coordinates
$(a, \alpha )$ on $\Q(1)$ as follows. If $(g, \gamma )$ 
are coordinates on $\GQ(1)$, then
$$
\begin{matrix}
\rho ^{*}(\alpha )&= (g+\gamma )^{-1}(a+\alpha )
(g+\gamma )]_{\bar 1}=\alpha \cr
\rho ^{*}(a)&=[(g+\gamma )^{-1}(a+\alpha )(g+\gamma )]_{\bar 0}=
a+2g^{-1}\gamma \alpha. \cr
\end{matrix}
$$
Notice that $\alpha $ is the function $\qtr $ on $\Q(1)$. The
following statement is obvious:

\begin{Theorem} The set of functions on $\Q(1)$ invariant
with respect to the adjoint action of $\GQ(1)$ coincides
with the set of functions of the form $\alpha\cdot f(a)+c$, where $f$ is
an arbitrary function on $\Cee $ and $c\in \Cee $.
\end{Theorem}

Therefore, any invariant function on $\Q(1)$ can be
expressed in terms of one invariant function $\alpha =\qtr  A$ and
one noninvariant function, $a$. The latter is, so to say, an
invariant of {\it second class}: $a$ is invariant on the
subsupermanifold singled out by the equation $\alpha =0$; besides, the
map under which $a$ passes to $\cpr  (a)$ --- a function on $\Cee =Q
(1)_{rd}$ --- is $\GQ(1)$-invariant.

\ssec {2.2. Invariant functions on $\Odd (1)$} Let 
$\widetilde{\Odd} (1)$ be the dense open subsupermanifold in $\Odd (1)$ singled
out by the equation
$$
\cpr  (a_{12}a_{21})\neq 0, 
$$
where $a_{ij}$ are standard coordinates on $\Odd (1)$. Let us
identify $\Cee ^{1|1}$ with the closed subsupermanifold of 
$\Odd (1)$
consisting of matrices of the form $\begin{pmatrix} \alpha & a \cr 1 & 0
\cr\end{pmatrix}$. Define the following 
even and an odd functions on $\Odd (1)$:
$$
g=a_{12}a_{21}-a_{11}a_{22}, \; \; \gamma =a_{11}+a_{22}.
$$
Thus, we have defined a map
$$
\pi :\widetilde{\Odd}(1)\tto \Cee ^{1|1}
\hookrightarrow \Odd (1)\; \; \text{ for which}\; 
\pi ^{*}(a)=g, \pi ^{*}(\alpha )=\gamma .
$$

\begin{Lemma} Any family of matrices $\varphi :{\cal U}\tto 
\widetilde{\Odd} (1)$ is
equivalent {\em (with respect to the $\GL (1|1)$-action)} to a family
$\pi \circ \varphi : {\cal U}\tto \widetilde{\Odd} (1)$. 
The set of all families equivalent to $\varphi $ is
mapped by $\pi $ into the set of all families $\begin{pmatrix} \alpha & a+
\varepsilon \alpha \cr 1 & 0 \end{pmatrix}$, where
$\varepsilon $ is an arbitrary odd function on ${\cal U}$.
\end{Lemma}

\begin{proof}
$$
\begin{pmatrix} 1 & -a_{22} \cr 0 & a_{21} \end{pmatrix}^{-1}
\begin{pmatrix} a_{11} & a_{12} \cr 
a_{21} & a_{22} \end{pmatrix}\begin{pmatrix} 1 & -a_{22} \cr 0 & a_{21}
\end{pmatrix} =\begin{pmatrix} a_{11}+a_{22} & a_{12}a_{21}-a_{11}a_{22} \cr 
1 &0\end{pmatrix}.
$$
If
$$
\begin{pmatrix} x & y \cr z & t \end{pmatrix}\begin{pmatrix} \alpha & a \cr 
1 &0 \end{pmatrix}=
\begin{pmatrix} \beta & b \cr 1 & 0 \end{pmatrix}\begin{pmatrix} x & y \cr 
z & t\end{pmatrix}
$$
where $x, t$ are even and $y, z$ are odd and $xt$ is invertible, 
then $\alpha =\beta$ since the supertrace is invariant and it is easy
to verify that $b-a=(2t^{-1}za)\alpha $.

On the other hand
$$
\begin{pmatrix} 1+a^{-1}\alpha \varepsilon & \varepsilon \cr 
a^{-1}\varepsilon & 1 \end{pmatrix}
\begin{pmatrix} \alpha & a \cr 1 & 0 \end{pmatrix}
\begin{pmatrix} 1+a^{-1}\alpha \varepsilon & \varepsilon \cr 
a^{-1}\varepsilon & 1 \end{pmatrix}^{-1}=
\begin{pmatrix} \alpha & a+2\varepsilon \alpha \cr 1 & 0 \end{pmatrix}.
$$
\end{proof}

\begin{Theorem} The set of $\GQ(1)$-invariant functions on
$\widetilde{\Odd} (1)$ coincides with the set of functions of the form
$\gamma \cdot h(g)+c$, where $h$ is an arbitrary function on
$\Cee \setminus \{0\}, c\in \Cee $.
\end{Theorem}

{\bf Proof}. Let $f$ be an invariant function. Denote $\pi ^{*}f$ by
$f'(a, \alpha )$. Then Lemma implies that $f'(a, \alpha )=f'(a+\varepsilon 
\alpha , \alpha )$
wherefrom
$$
f'(a, \alpha )=\alpha h(a)+c,\quad \quad  f=f'\circ \pi ^{*}=\gamma h(g)+c.
$$
Conversely, if $f=\gamma h(g)$, then by Lemma $\gamma $ is an
invariant and the noninvariance of $g$ is equivalent to
the replacement of $g$ by $g+\delta \gamma $ wherefrom
$$
\gamma h(g+\varepsilon \gamma )=\gamma h(g)+\gamma \varepsilon \gamma h'(g)=
\gamma h(g).
$$

{\bf Example}. It is easy to see that $\str A^{2n}=0, \str 
A^{2n+1}=(2n+1)\gamma g^{n}$.

\ssec{2.3} As follows from Lemma 2.2, the matrix 
$\begin{pmatrix} a_{11} & a_{12} \cr a_{21} & a_{22} \end{pmatrix}$ is
equivalent to any of the matrices of the form
$\begin{pmatrix} \gamma & g+\varepsilon \gamma \cr 1 & 0 \end{pmatrix}$ so for 
description of invariant functions on
$\Odd (1)$ instead of $g$ we could have used any function
$g'=g+\varepsilon \gamma $.

The arbitrariness in the choice of $g$ is connected with a slightly more 
general
question: what data on ${\cal U}$-family of matrices $A$ (in other words, on a
matrix-valued function
$A$ on ${\cal U}$) should be given in order to enable us to compute the values
of all the invariant functions. Let $a$ be a function on ${\cal U}$ such that
the value of any invariant function $f=\gamma h(g)$ on $A$ is equal to $\str A
\cdot h(a)$.  Then, in particular, $\str A^{3}=3a \cdot \str A$. We get a 
condition on $a$ which is equivalent to the equation
$$
[a-g(A)]\str A=0.
$$
It turns out that fulfilment of this equation suffices.

\begin{Theorem} Let $A$ be a matrix-valued function on ${\cal U}$ with
values in $\widetilde{\Odd} (1)$ and $a$ an even function 
on ${\cal U}$ such that
$$
\str A^{3}=3a \cdot \str A.
$$
Then for any function $h$ on $\Cee \setminus \{0\}$ we have
$$
\gamma (A)h(g(A))=\str A \cdot h(a).
$$
\end{Theorem}

\begin{proof} Let $u_{1}, \dots , u_{k}$ be even and $v_{1}, \dots , v_{l}$ odd
local coordinates on ${\cal U}$. The condition $[a-g(A)]\str A=0$ means that 
for
any fixed value of coordinates $u_{1}, \dots , u_{k}$ either
$\str A=0$ or $\Delta =a-g(A)$ belongs to the algebra generated by
$v_{1}, \dots , v_{l}$.

In the first case both sides vanish; in the
second case we have
$$
\gamma (A) \cdot h(a)=\gamma (A)\left\{h(g(A))+\sum^{l}_{i=1}
\frac{\Delta ^{i} }{ i!}\frac{d^{i}b}{dz^{i}}(g(A))\right\}=
\gamma (A) \cdot h(g(A)).
$$
\end{proof}

\begin{Corollary} If $A_{1}$ and $A_{2}$ be two ${\cal U}$-families of
matrices from $\widetilde{\Odd} (1)$ then they are indistinguishable by
$\GL (1|1)$-invariant functions on $\widetilde{\Odd} (1)$ 
if and only if
$\str A_{1}=\str A_{2}$ and $\str A^{3}_{1}=\str 
A^{3}_{2}$.
\end{Corollary}

\begin{Remark} 1) The condition $[a-g(A)]\str A=0$ does not
generally imply the equality $a=g(A)+\varepsilon (A)\gamma (A)$ even if
$\gamma (A)\neq 0$; for example, take $A=\begin{pmatrix} v_{1}v_{2}v_{3} & 1 
\cr
1 & 0 \end{pmatrix}$ and $a=1+v_{1}v_{2}$, 
where $v_{1}, v_{2}, v_{3}$ are odd coordinates on ${\cal U}$. Nevertheless, if
$g'$ is an odd function (or family of functions) on an open
subsupermanifold $U\subset \widetilde{\Odd} (1)$ for which 
$\gamma \cdot g'=\gamma \cdot g|_{U}$, then
$g'=g+\varepsilon \gamma $ for an odd function (or family of functions) 
$\varepsilon $ on $U$.
To make sure of this it suffices to take $\gamma $ and $a_{11}-a_{22}$
for odd coordinates on $\Odd (1)$.

2) It is not difficult to verify that everything
said in this subsection can be translated almost
literally to the case of $\Q(1)$ with inessential
distinctions: first, there is a distinguished
even function, $a$, on $\Q(1)$, and second, the appearance of a nonzero
summand $\varepsilon a$ is only possible in the presence of odd
parameters.
\end{Remark}

Thus, both on $\Q(1)$ and on $\Odd (1)$ there exists a
pair consisting of one odd invariant function $\tau\; (\alpha $ or $\gamma )$
and one even noninvariant function $t\; (a$ or $g)$ such that
any invariant function can be expressed in the form
$$
\tau h(t)+{\rm const}, \; \text{
where $\tau$ and $t$ are polynomials}.
$$

In this section the study of invariant functions on
$\Q(1)$ and $\Odd (1)$ was reduced to the study of functions
invariant with respect to the action of the supercommutative
supergroup $\Cee ^{0|1}$ on $\Cee ^{1|1}$.

The next section is devoted to the $S_{n}\vdash \Cee ^{0|n}$-invariant
functions on $\Cee ^{n|n}$. These functions play a similar role in the study
of $\GQ(n)$-invariant functions on $\Q(n)$ and $\GL
(n|n)$-invariant functions on $\Odd (1)$.

\section* {\protect \S 3. Invariant functions on $\tilde{\Cee }^{n|n}$}

\ssec{3.1} The discrete group of permutations $S_{n}$ acts on
the supermanifold $\Cee ^{n|n}$ with coordinates
$a_{1}, \dots , a_{n}, \alpha _{1}, \dots , \alpha _{n}$ by permuting indices 
of
the coordinate functions of the same parity and the supercommutative supergroup
$\Cee ^{0|n}$ with coordinates $\varepsilon _{1}, \dots , \varepsilon _{n}$ 
acts
as follows. The action of $\Cee ^{0|n}$ is the morphism
$\rho : \Cee ^{0|n}\times \Cee ^{n|n}\rightarrow 
\Cee ^{n|n}$ given by the formulas
$$
\rho ^{*}(\alpha _{i})=\alpha _{i}, \rho (a_{i})=a_{i}+\varepsilon _{i}
\alpha _{i}, \; \; i=1, \dots
, n.
$$
Thus, the semidirect product of these supergroups, the supergroup $S_{n}\vdash
\Cee ^{0|n}$,  acts on $\Cee ^{n|n}$.  Denote by $\tilde{\Cee }^{n}$ the open
subset of $\Cee ^{n}$ consisting of $n$-tuples of pair-wise distinct complex
numbers and by $\tilde{\Cee }^{n|n}$ the  subsupermanifold of $\Cee ^{n|n}$
whose underlying is ??. In what follows  we will call the $S_{n} \vdash \Cee
^{0|n}$-invariant functions on  $\tilde{\Cee }^{n|n}$ just {\it invariants} and
$S_{n}$-invariant functions {\it symmetric functions} (inside of this section 
we will not encounter other types
of invariance).

\begin{Theorem} The set of invariant functions coincides
with the set of functions of the form
$$
f=f_{0}\sum^{}_{i}\alpha _{i}f_{1}(a_{i})+\sum^{}_{i<j}
\alpha _{i}\alpha _{j}f_{2}(a_{i}, a_{j})+\dots
+\alpha _{1}\dots
\alpha _{n}f_{n}(a_{1}, \dots
, a_{n}), 
\eqno{(3.1)}
$$
 where $f_{k}$ is an even function in $k$ even variables 
skew-symmetric with respect to these variables for $k>1$.
\end{Theorem}

\begin{proof} Let us express an arbitrary invariant
function in the form
$$
\begin{matrix}
f&=g_{0}(a_{1}, \dots
, a_{n})+\sum \alpha _{i}g_{i}(a_{1}, \dots
, a_{n})+\dots
+ \cr
&+\sum^{}_{i<j}\alpha _{i}\alpha _{j}g_{ij}(a_{1}, \dots
, a_{n})+\dots
+\alpha _{1}\dots
\alpha _{n}g_{12\dots
n}(a_{1}, \dots
, a_{n}). \end{matrix}
$$
 The invariance with respect to $\Cee ^{0|n}$ 
is equivalent to the
conditions
$$
\alpha _{i} \cdot \pder{a_{i}}f=0 \text{ for all $i$}
$$
 wherefrom we see that $g_{i_{1}\dots i_{s}}$ only depend on the
$a_{i_{1}}\dots a_{i_{s}}$. If $\delta \in S_{n}$ is such that $\delta (i_{1},
\dots , i_{s})=(1, \dots , s)$ then comparing the coefficient of $\alpha
_{1}\dots \alpha _{s}$ for $f$ and $\delta f$ we get $g_{i_{1}\dots
i_{s}}=g_{1\dots s}$ and if $\delta $ preserves $s+1, \dots , n$, then
$$
\begin{matrix}
&\alpha _{1}\dots
\alpha _{s}g_{1\dots
s}(a_{1}\dots
a_{s})=\alpha _{\delta _{1}}\dots
\alpha _{\delta _{s}}g_{1\dots
s}(a_{s_{1}}\dots
a_{s_{n}})= \cr
&=(-1)^{p(\delta )}\alpha _{1}\dots
\alpha _{s}g_{1\dots
s}(a_{s_{1}}\dots
a_{s_{n}}), \end{matrix} 
$$
where $p(\delta )$ is the parity of the permutation 
$\delta $ which is
equivalent to skew symmetricity of $f_{s}=g_{1\dots
s}$. The invariance of functions of the form (3.1) is obvious.
\end{proof}

\begin{Remark} Clearly, $f$ is a rational or polynomial function if and only if
all the functions $f_{1}, \dots , f_{n}$ are of the same class (rational or
polynomial). Set $\tau_{k}=\sum \alpha _{i}a^{k-1}_{i},  t_{k}=\sum a^{k}_{i}$.
Theorem implies that functions $\tau_{k}$ are invariant whereas functions
$t_{k}$ are not.
\end{Remark} 

It is not difficult to see that the algebra of invariant functions does not 
have
a finite set of generators regardless of the class of functions (polynomials,
rational or analytic functions, $\dots )$ with which we work: it follows from
(3.1) that any invariant function differs by a constant from a nilpotent 
function
--- an element of the ideal generated by functions $\alpha _{1}, \dots , \alpha
_{n}$ and, therefore, a finite set of invariant functions only generates a
finite dimensional subspace even if we admit rational expressions, in infinite
dimensional (since $\tau_{k}$ are linearly independent) space of invariant
functions.

Let us prove a superanalog of the main theorem on
symmetric functions.

\ssbegin{3.2}{Theorem} If $f(a, \alpha )$ is a symmetric function, then

$1)$ there exists a unique function $g(x_{1}, \dots , x_{n}, \xi _{1}, \dots ,
\xi _{n})$ defined on an open subsupermanifold of $\Cee ^{n|n}$ such that
$g(t_{1}, \dots , t_{n}, \tau_{1}, \dots , \tau_{n})=f$; 

$2)$ if $f$ is a polynomial
function, then so is $g$ and if $f$ is a rational function
then so is $g$.
\end{Theorem}

We will need an auxiliary statement. Set
$$
M(a_{1}, \dots
, a_{n})=\begin{pmatrix} 1 & \dots & 1 \cr 
\vdots&\vdots&\vdots \cr 
a^{n-1}_{1} & 
\dots & a^{n-1}_{n} \end{pmatrix}.
$$

\begin{Lemma}
$$
M^{-1}=\begin{pmatrix} \prod\limits _{i\neq 1}(a_{1}-a_{i})^{-1} & 0 \cr 0 & 
\prod\limits _{i\neq n}(a_{n}-a_{i})^{-1}\end{pmatrix} M', 
\; \; \text{{\it where $M'$ is a polynomial matrix}}.
$$
\end{Lemma}

{\bf Proof of Lemma}.  Let $M'$ be a matrix whose $s$-th row consists 
of coefficients of the polynomial $\prod _{1\le i\le n}(x-a_{i})$ 
written in order of increase of the power of $x$.  Then
$$
(M'M)_{kl}=\prod \limits _{1\le i\le n, \; \; i\neq k}(a_{l}-a_{i}). 
\qed
$$

{\bf Proof of Theorem}.  Let us express $\alpha _{1}, \dots , \alpha 
_{n}$ through $\tau_{1}, \dots , \tau_{n}$ and $a_{1}, \dots , a_{n}$:
$$
\begin{pmatrix} \alpha _{1} \cr \vdots \cr \alpha _{n} \end{pmatrix}
=M^{-1}(a_{1}, \dots , a_{n})\begin{pmatrix} \tau_{1} \cr \vdots \cr \tau_{n}
\end{pmatrix}
\eqno{(3.2)}
$$
Formula (3.2) implies that $\tau_{1}, \dots , \tau_{n}, a_{1}, \dots , 
a_{n}$ is a global coordinate system on $\tilde{\Cee }^{n|n}$.  Let 
$f(a, \alpha )$ be a symmetric function; let us express it in terms of 
$\tau, a$:
$$
f=c_{0}(a_{1}, \dots , a_{n})+\sum^{}_{1\le s\le n}\; \; \sum^{}_{1\le
i_{1}<\dots <i_{s}\le n}\tau_{i_{1}}\dots
\tau_{i_{s}}c_{i_{1}\dots i_{s}}(a_{1}, \dots , a_{n})\eqno{(3.3)}
$$
Since the $\tau_{i}$ are $S_{n}$-invariant, it follows that all the 
coefficients $c_{i_{1}\dots i_{s}}$ are symmetric in $a_{1}, \dots , 
a_{n}$ and therefore can be expressed in terms of $t_{1}, \dots , 
t_{n}$ because the Jacobian of the map
$$
(a_{1}, \dots
, a_{n})\mapsto (\sum a_{i}, \sum a^{2}_{i}, \dots
, \sum a^{n}_{i})
$$
is an invertible on $\tilde{\Cee }^{n}$ matrix $M(a_{1}, \dots , 
a_{n})$ and determines a diffeomorphism of $\tilde{\Cee }^{n}/S_{n}$ 
with an open submanifold $U\subset \Cee ^{n}$ (the complement to the 
set of zeros of a polynomial).

Let us establish that a symmetric function $f(a, \alpha )$ can be expressed in
the form $g(t, \tau)$ in a unique way. The functions $a_{1}, \dots , a_{n},
\tau_{1}, \dots , \tau_{n}$ constitute a global coordinate system on
$\tilde{\Cee }^{n|n}$ and, therefore, if
$$
d_{0}(t_{1}, \dots , t_{n})+\sum^{}_{1\le s\le n}\sum^{}_{1\le i_{1}<\dots
<i_{s}\le n}\tau_{i_{1}}\dots
\tau_{i_{s}}d_{i_{1}\dots i_{s}}(t_{1}, \dots
, t_{n})=0,
$$
then $d_{0}$ and all $d_{i_{1}\dots i_{s}}$ vanish on $U$.

If $f(a, \alpha )$ is a rational function, then $c_{i_{1}\dots i_{s}}$ 
in (3.3) are also rational since they are linear combinations of 
rational functions $f_{s}(a_{i_{1}}, \dots , a_{i_{s}})$ with rational 
coefficients --- polynomials in matrix elements of $M^{-1}$ and then, 
as immediately follows from theorem on symmetric polynomials, the 
functions $d_{i_{1}\dots i_{s}}$ determined from the condition 
$c_{i_{1}\dots i_{s}}(a_{1}, \dots , a_{n})=d_{i_{1}\dots 
i_{s}}(t_{1}, \dots , t_{n})$ are also rational.

Let now
$$
f(a, \alpha )=f_{0}+\sum^{}_{1\le s\le n}\; \; \sum^{}_{1\le i_{1}<\dots
<i_{s}\le n}\alpha _{i_{1}}\dots
\alpha _{i_{s}}f_{i_{1}\dots
i_{s}}
$$
where $f_{0}, f_{i_{1}\dots
i_{s}}$ are polynomials in $a_{1}, \dots
, a_{n}$. Then
$$
f_{i_{1}\dots
i_{s}}(a_{1}, \dots
, a_{n})=q_{i_{1}\dots
i_{s}}(a_{1}, \dots
, a_{n})\prod _{1\le k<l\le s}(a_{i_{k}}-a_{i_{l}})
\text{ for $s>1$}, 
$$
where $q_{i_{1}\dots i_{s}}$ are also polynomials since $f_{i_{1}\dots 
i_{s}}$ are skew-symmetric with respect to permutations of $a_{i_{1}}, 
\dots , a_{i_{s}}$.

By Lemma
$$
\alpha _{j}=\prod _{1\le s\le n}(a_{j}-a_{s})^{-1}
\sum^{}_{k}(M')_{jk}\tau_{k};
$$
hence, 
$$
\begin{matrix}
&\alpha _{i_{1}}\cdot \dots
\cdot \alpha _{i_{s}}\prod _{1\le k<l\le s}(a_{i_{k}}-a_{i_{l}})= \cr
&=\sum^{}_{j_{1}, \dots
, j_{s}}\left(\prod _{i\le k<l\le s}(a_{i_{k}}-a_{i_{l}})\right)
\left(\prod _{1\le r\le n, r\neq i_{1}}(a_{i_{1}}-a_{r})^{-1}\right)\dots
\left(\prod _{1\le r\le n, r\neq i_{s}}(a_{i_{s}}-a_{r})^{-1}\right)\cdot \cr
& (M')_{i_{1}j_{1}}\dots
(M')_{i_{s}j_{s}}\tau_{j_{1}}\dots
\tau_{j_{s}}=\prod _{i\le k<l\le n}(a_{k}-a_{l})^{-1}\cdot P_{i_{1}\dots
i_{s}}, \end{matrix}
$$
where $P_{i_{1}\dots i_{s}}$ is a polynomial in $a_{1}, \dots , a_{n}, 
\tau_{1},
\dots , \tau_{n}$.

Indeed, the factor $a_{i_{k}}-a_{i_{l}}$ appears in
this expression three times: in $\prod (a_{i_{k}}-a_{i_{l}})$, in 
$\prod (a_{i_{s}}-a_{r})^{-1}$ and in
$\prod (a_{i_{l}}-a_{r})^{-1}$ and, therefore, its total power is equal to 
$-1$.
This implies that
$$
f(a, \alpha )=\prod _{k<l}(a_{k}-a_{l})^{-1}\cdot P(a, \tau), 
$$
 where $P$ is a polynomial. Since $f$ is a symmetric function, 
$$
P(a, \tau)=f(a, \alpha )\cdot \prod \limits_{k<l}(a_{k}-a_{l})
$$
is skew-symmetric and, therefore, is divisible by $a_{k}-a_{l}$, i.e., 
$P=(\prod \limits_{k<l}(a_{k}-a_{l}))\Q (a, \tau)$, where $\Q(a, 
\tau)$ is a polynomial in $a_{1}, \dots , a_{n}, \tau_{1}, \dots , 
\tau_{n}$ symmetric with respect to $a_{1}, \dots , a_{n}$.  Then 
$f(a, \alpha )=\Q(a, \tau)$ and it remains to express $\Q(a, \tau)$ in 
the form of a polynomial in $t_{1}, \dots , t_{n}, \tau_{1}, \dots , 
\tau_{n}$ which is possible thanks to a theorem on symmetric 
polynomials.

\begin{Remark} The theorem proved above will look quite naturally if 
we allow functions in non-homogeneous argument assuming that there 
exists a unique decomposition of a non-homogeneous (with respect to 
parity) $x$ into the sum of an even and an odd summands.  Then the 
theorem means precisely the following:

In the algebra of all (polynomial, rational, etc.)
functions in non-homogeneous arguments
$x_{1}=a_{1}+\alpha _{1}, \dots
, x_{n}=a_{n}+\alpha _{n}$ the subalgebra of symmetric
functions is generated by
$$
y_{1}=\sum (a_{i}+\alpha _{i}), y_{2}=\sum (a^{2}_{i}+
2a_{i}\alpha _{i}), \dots
, y_{n}=\sum (a^{n}_{i}+n\alpha _{i}a^{n-1}_{i}).
$$
\end{Remark}

\ssec{3.3} Let us write the conditions that single out the
invariant functions from the set of all functions in
$t_{1}, \dots
, t_{n}, \tau_{1}, \dots
, \tau_{n}$, i.e., from all the symmetric
functions on $\tilde{\Cee }^{n|n}$.

Let us denote by $\check \Cee ^{n|n}$ the image of 
$\tilde \Cee ^{n|n}$ under the
map to $\Cee ^{n|n}$ given by functions $t_{1}, \dots
, t_{n}, \tau_{1}, \dots
, \tau_{n}$. Any function $h$ on $\check \Cee ^{n|n}$ such that
$h(t_{1}(a, \alpha ), \dots
, \tau_{n}(a, \alpha ))$ is an invariant function will be referred to a
{\it balanced function}.

\begin{Lemma}  A function $h(u, \xi )$ is balanced if and only if
$$
\sum^{n}_{s=1}s\cdot \tau_{i+s-1}\pderf{h}{u_{s}}(t_{1}, \dots
, t_{n}, \tau_{1}, \dots
, \tau_{n})=0 \; for\; i=1, \dots
, n. 
$$
\end{Lemma}

\begin{proof} Let us transform the conditions of
$\Cee ^{0|n}$-invariance
$$
\alpha _{i}\pderf{f}{a_{i}}=0, \; \; i=1, \dots
, n
$$
 with the help of an invertible matrix $M(a_{1}, \dots
, a_{n})$. This produces equivalent but symmetric conditions
$$
\begin{matrix}
&\sum^{}_{j}\alpha _{j}a^{i-1}_{j}\pderf{f}{a_{j}}=0 \; \; 
\text{for $i=1, \dots
, n$.} \cr
&\sum^{}_{j}\alpha _{j}a^{i-1}_{j}\pder{a_{i}}h(t_{1}, 
\dots
, t_{n}, \tau_{1}, \dots
, \tau_{n})= \cr
&=\sum^{}_{s}\sum^{}_{j}s\alpha _{j}a^{i+s-2}_{j}\pderf{h}{u_{s}}(t_{1}, \dots
, t_{n}, \tau_{1}, \dots
, \tau_{n})= \cr
=&=\sum^{}_{s}s\tau_{i+s-1}\pderf{h}{u_{s}}(t, \tau). \end{matrix}
$$
\end{proof}

Lemma gives us a possibility to describe polynomial
balanced functions, i.e., invariant polynomials on $\Cee ^{n|n}$.

Actually, invariant polynomials are already described
by Sergeev in \cite{S}. Here we give a new proof based on two
ideas: 

1) the passage to symmetric polynomials in an
infinite number of indeterminates makes variables
independent (cf. \cite{M}) and

2) homology of the vector field $\sum \xi _{i}\pder{u_{i}}$ on 
$\Cee ^{n|n}$
consist of constants (we can
speak about homology because $(\sum \xi _{i}\pder{u_{i}})^{2}f=0$).

\begin{Theorem} The algebra of invariant polynomials on
$\Cee ^{n|n}$ is generated by functions $\tau_{1}, \dots
, \tau_{k}, \dots$: any
invariant polynomial can be uniquely expressed in the form
$$
\sum^{}_{0\le s\le n}\; \; \sum^{}_{i_{1}<\dots
<i_{s}}c_{i_{1}, \dots
, i_{s}}\tau_{i_{1}}\dots
\tau_{i_{s}}
$$
where only a finite number of coefficients is nonzero.  All the 
relations between $\tau_{1}, \dots , \tau_{n}, \dots $ are generated 
by identities of the form
$$
\tau_{i_{1}}\dots
\tau_{i_{n+1}}=0.
$$
\end{Theorem} 

\begin{proof} Denote by ${\cal P}_{n}$ the graded superalgebra of
symmetric polynomials in $a_{1}, \dots
, a_{n}$, $\alpha _{1}, \dots
, \alpha _{n}$, where $\deg
\alpha _{i}=\deg a_{i}=1$, and by $\pi _{n}: {\cal P}_{n+1}\tto 
{\cal P}_{n}$ the projection that sends
$a_{n+1}$ and $\alpha _{n+1}$ to zero and the remaining generators to
their namesakes; denote by ${\cal P}_{\infty }$ the projective limit of 
${\cal P}_{n}$ in the
category of $\Zee $-graded rings. Theorem 3.2 implies that the
$\Zee $-graded supercommutative superalgebra ${\cal P}_{\infty }$ 
is isomorphic
to the graded superalgebra $Q$ of polynomials in two
countable sets of generators: the even ones, $f_{1}, \dots
, f_{k}, \dots
$
and the odd ones, $\varphi _{1}, .., \varphi _{k}, \dots
$, where $\deg f_{k}=\deg \varphi _{k}=k$.

Indeed, $t_{1}, \dots , t_{n}, \tau_{1}, \dots , \tau_{n}$ are functionally,
hence, algebraically, independent in ${\cal P}_{n}$ and, therefore, $\{t_{i}, 
\tau_{i}\;  i\in\Zee \}$ are algebraically independent in ${\cal P}_{\infty
}$.

On $Q$, the action of differentiations
$$
v_{k}=\sum\limits^{\infty }_{s=1}s\varphi _{k+s}\pder{f_{s}}
$$ 
is 
well-defined since every element of
$Q$ only contains a finite number of generators $f_{s}$. Let us
consider diagrams
$$
\begin{matrix}\begin{matrix}
{\cal P}_{n+1} &\buildrel \sum\limits_i \alpha
_{i}a^{k-1}_{i}\pder{a_{i}}\over{\tto}&  {\cal P}_{n+1}\cr
\downarrow&&\downarrow\cr
 {\cal P}_{n}&\buildrel \sum\limits_i\alpha _{i}a^{k-1}_{i}
\pder{a_{i}}\over{\tto}& {\cal P}_{n}\end{matrix}
&
\begin{matrix}
{\cal P}_{n}&\buildrel \sum\limits_i\alpha
_{i}a^{k-1}_{i}\pder{a_{i}}\over{\tto}&  {\cal P}_{n}\cr
\downarrow&&\downarrow\cr 
Q &\buildrel v_{k}\over{\tto}& Q\end{matrix}\end{matrix}
$$
 The commutativity of the first one is obvious and the
commutativity of the second one is proved together with
Lemma 3.3; therefore, the limit of the subsuperalgebras
of invariant polynomials in ${\cal P}_{n}$ is a subsuperalgebra of 
${\cal P}_{\infty }$
which under the isomorphism ${\cal P}_{\infty }\cong Q$ turns into the
subsuperalgebra singled out by equations $v_{k}f=0$.

Let $R\in {\cal P}_{n}$ be an invariant polynomial of degree $k$; by
Theorem 3.1 it is determined by a set $r_{0}, \dots
, r_{n}$, where $r_{i}$
is a polynomial in $i$ variables and, therefore, for any $m>n$
there exists an invariant polynomial $R_{m}\in {\cal P}_{n}$ given by the
same functions $r_{0}, \dots
, r_{k}$. The images of $R_{m}$ in $Q$ stabilize
for $m>k$ giving rise to a polynomial $S$. The condition
$$
0=v_{k+1}S=\sum^{}_{1\le j\le k}j\varphi _{j+k+1}\pderf{S}{f_{j}}
$$
 implies $\pderf{S}{f_{k}}=0$ since $S$ 
does not depend on $\varphi _{j}, f_{j}$, where
$j>k$, and the generators of $Q$ are independent. Therefore, $S$
is a polynomial in $\varphi _{1}, \dots
, \varphi _{k}$, i.e., 
$$
S=c_{0}+\sum^{}_{j_{1}<j_{2}<\dots
<j_{i}, j_{1}+\dots
+j_{i}\le k}c_{j_{1}\dots
j_{i}}\varphi _{j_{1}}\dots
\varphi _{j_{i}}
$$
 and
$$
R=S(\tau_{1}, \dots
, \tau_{k})=c_{0}+\sum^{n}_{i=1}\sum^{}_{j_{1}<j_{2}<\dots
<j_{i}, j_{1}+\dots
+j_{i}\le k}c_{j_{1}\dots
j_{i}}\tau_{j_{1}}\dots
\tau_{j_{i}}, 
$$
where we have taken into account that $\tau_{i_{1}}\dots
\tau_{i_{n+1}}=0$ in
${\cal P}_{n}$.

It remains to prove that in ${\cal P}_{n}$ all the polynomials
$\tau_{i_{1}}\dots
\tau_{i_{s}}$ with $i_{1}<i_{2}<\dots <i_{s}$ and $s\le n$ are linearly
independent. In ${\cal P}_{n}$, introduce a grading setting $\deg 
\alpha _{i}=\deg a_{i}=i$. Then the lowest term of $\tau_{i_{1}}\dots
\tau_{i_{s}}$ is equal to
$\alpha _{1}\dots
\alpha _{s}a^{i_{s}}_{1}a^{i_{s-1}}_{2}\dots a^{i_{1}}_{s}$. The equality
$\tau_{i_{1}}\dots
\tau_{i_{n+1}}=0$ is obvious. \end{proof}

\ssec{3.4} By analogy with the case $n=1$ considered in subsec. 2.3 we see that
the action of $\Cee ^{0|n}$ on the set of functions $t_{i}$ produces a new set
of even functions equally suitable for calculation of invariant functions. To
study this nonuniqueness it is useful to replace the collection
$t_{1}, \dots , t_{n}$ with another collection of generators in the algebra of
symmetric functions in $a_{1}, \dots , a_{n}$. Set
$$
s_{j}=(-1)^{j-1}\sum^{}_{1\le i_{1}<\dots
<i_{j}\le n}a_{i_{1}}\dots
a_{i_{j}}.
$$
Then the sets $\{t_{i}\}$ and $\{s_{j}\}$ are  expressed in terms of each other
polynomially and the functions
$(s_{1}, \dots , s_{n}, \tau_{1}, \dots , \tau_{n})$ constitute a coordinate
system on
$\check \Cee ^{n|n}$.

\begin{Lemma}
$$
\tau_{n+k}=\sum^{}_{1\le j\le n}\tau_{n+k-j}s_{j}
\eqno{(3.4)}
$$
\end{Lemma}

\begin{proof} The equality
$$
\prod _{1\le i\le n}(x-a_{i})=x^{n}-\sum^{}_{1\le j\le n}s_{j}x^{n-j}
$$
implies $a^{n}_{i}=\sum^{}_{1\le j\le n}s_{j}a^{n-j}_{i}$, 
i.e., the coefficients of $\alpha _{i}$ in
the left- and right-hand sides of (3.4) coincide.

In what follows an arbitrary supermanifold of
parameters is denoted by ${\cal U}$ and the sheaf of ideals in
${\cal F}({\cal U})$ generated by odd functions is denoted by 
$I_{{\cal U}}$ or just
by $I$.

Two families of ${\cal U}$-points
$$
\tilde \Cee ^{n|n}: \varphi _{1}: {\cal U}\tto 
\tilde \Cee ^{n|n}, \; \; i=1, 2
$$
will be called {\it equivalent} if
$$
\psi ^{*}_{1}(\alpha _{i})=\psi ^{*}_{2}(\alpha _{i}), 
\psi ^{*}(\alpha _{i})[\psi ^{*}_{1}(a_{i})-\psi ^{*}_{2}(a_{i})]=0 
\; \; \text{for all $i$.}
$$
Theorem 3.1 implies that the values of any invariant
function $f$ on equivalent families $\psi _{1}$ and 
$\psi _{2}$ coincide, i.e., 
$\psi ^{*}_{1}(f)=\psi ^{*}_{2}(f)$.
\end{proof}

\begin{Theorem} Let $\varphi : {\cal U}\tto 
\tilde \Cee ^{n|n}$ be a family of ${\cal U}$-points of
$\tilde \Cee ^{n|n}$ and $g_{1}, \dots
, g_{n}$ functions on ${\cal U}$ such that
$$
\varphi ^{*}(\tau_{n+k})=\sum^{n}_{i=1}g_{i}\varphi ^{*}
(\tau_{n+k-i})\; for\; k+1, \dots
, n 
\eqno{(3.5)}
$$
 and either
$$
\varphi ^{*}(\alpha _{i})\neq 0\; for\; all\; i\le n 
$$
 or
$$
g_{j}=\varphi ^{*}(s_{j}) (\mod  I_{{\cal U}})\; for\; all\; j\le n.
$$
 Then there exists a unique equivalent to $\varphi $ family of
morphisms $\psi : {\cal U}\tto \tilde \Cee ^{n|n}$ 
such that $g_{j}=\psi ^{*}(s_{j})$ for all $j$.
\end{Theorem}

\begin{proof} First, consider an homogeneous system of
equations
$$
\sum^{n}_{j=1}\Delta _{j}(\sum^{n}_{i=1}\beta _{i}b^{n+k-j-1}_{i})=0 \; \; 
\text{for $k=1, \dots
, n $}
\eqno{(3.6)}
$$
 in which $\beta _{i}$ are odd functions on ${\cal U}$ and 
$\Delta _{j}$ and $b_{i}$ are even
functions and $\prod \limits_{i\neq j}(b_{i}-b_{j})$ is an 
invertible element in the
algebra of functions.

Denote by $M'(b_{1}, \dots
, b_{n})$ an $n\times n$-matrix whose $i$-th
column consists of coefficients of the polynomial
$\prod\limits _{i\neq j}(x-b_{j})$ written down in order of decreasing of the
power of $x$. Then similarly to Lemma 3.1 we have
$$
\begin{pmatrix} b^{n-1}_{1} & \dots & 1 \cr 
\vdots&\vdots&\vdots\cr 
b^{n-1}_{n} & 
\dots & 1 \end{pmatrix} M'(b_{1}, \dots
, b_{n})=\begin{pmatrix} \prod\limits _{i\neq 1}(b_{1}-b_{i}) & \dots & 0 \cr 
& \ddots&\cr
 0 & \dots & \prod \limits_{i\neq n}(b_{n}-b_{i}) \end{pmatrix}\eqno{(3.7)}
$$
In particular, both factors (matrices) are invertible.

Since
$$
0=\sum^{}_{j}\Delta _{j}(\sum^{}_{i}\beta _{i}b^{n+s-j-1}_{i})=
\sum^{}_{i}[\beta _{i}(\sum^{}_{j}\Delta _{j}b^{n-j-1}_{i})]b^{s}_{i}, 
$$
 then system (3.6) is equivalent to the system
$$
\beta _{i}(\sum^{}_{j}\Delta _{j}b^{n-j-1}_{i})=0\; \; i=1, \dots
, n
\eqno{(3.8)}
$$
Set $\Delta ^{0}_{j}=g_{j}-\varphi ^{*}(s_{j})$. Then $\Delta ^{0}_{1}, 
\dots
, \Delta ^{0}_{n}$ satisfies 
system (3.6) and, therefore, it satisfies (3.8) with
$\beta _{i}=\varphi ^{*}(\alpha _{i}), b_{i}=\varphi ^{*}(a_{i})$. 
If $\sum^{}_{j}\Delta _{j}b^{n-j-1}_{i}\not\equiv 0 ~(\mod  
I_{{\cal U}})$ for some
$i$, then (3.8) implies $\beta _{i}=0$ and, therefore, under conditions
of Theorem we always have $g_{j}\equiv \varphi ^{*}(s_{j}) (\mod  
I_{{\cal U}})$.

Let us show that for any $k\ge 0$ there exists a family
$\psi ^{(k)}$ equivalent to $\varphi $ such that $\psi ^{(k)*}(s_{j})
\equiv g_{j} (\mod  I^{k+1}_{{\cal U}})$.
For $k=0$ it suffices to take $\psi ^{(0)}=\varphi $. Let 
$\Delta ^{(k)}_{j}=\psi ^{(k)}(s_{j})-g_{j}$
and
$$
\begin{pmatrix} \Box  ^{(k)}_{1} \cr \vdots \cr \Box 
^{(k)}_{n}\end{pmatrix} = M^{-1}(\psi ^{(k)*}(a_{1}), \dots
, \psi ^{(k)*}(a_{n}))\begin{pmatrix} \Delta ^{(k)}_{1} \cr \vdots \cr 
\Delta ^{(k)}_{n} \end{pmatrix}.\eqno{(3.9)}
$$
 Then $\Box  ^{(k)}_{i}\in I^{k+1}_{{\cal U}}$ and 
$\Box  ^{(k)}_{i}\Box  ^{(k)}_{l}\equiv 0 (\mod  
I^{k+2}_{{\cal U}})$ and
therefore setting $\psi ^{(k+1)*}(a_{i})=\psi ^{(k)*}(a_{i})+
\Box  _{i}$ we get
$$
\begin{matrix}
\psi ^{(k+1)*}(s_{j})&\equiv \psi ^{(k)*}(s_{j})-\sum^{}_{i}
\psi ^{(k)*}\{s_{j-1}(a_{1}, \dots
, \hat a_{i}, \dots
, a_{n})\}\Box  {}^{(k)}\equiv \cr
&\equiv g_{j} (\mod  I^{k+2}_{{\cal U}}) \end{matrix}
$$
(here, as usual, hat over a symbol manifests its absence).

Let $\beta _{i}=\varphi ^{*}(\alpha _{i}), b_{i}=\varphi ^{(k)*}(a_{i})$, 
then (3.7) and (3.9)
imply
$$
\begin{pmatrix} b^{n-1}_{1} & \dots & 1 \cr 
\vdots&\vdots&\vdots \cr 
b^{n-1}_{n} & \dots & 1 \end{pmatrix}
\begin{pmatrix} \Delta ^{(k)}_{1} \cr \vdots \cr \Delta ^{(k)}_{n}
\end{pmatrix} =
\begin{pmatrix} \prod\limits _{i\neq 1}(b_{1}-b_{i})\Box  {}^{(k)}_{1} \cr 
\vdots \cr
\prod \limits_{i\neq n}(b_{n}-b_{i}) \Box  {}^{(k)}_{n} \end{pmatrix}.
$$
 Since $\psi ^{(k)}$ is equivalent to $\varphi $, then 
$\Delta ^{(k)}_{j}$ satisfies (3.6)
and (3.8) with $\beta _{i}=\varphi ^{*}(\alpha _{i}), b_{i}=
\varphi ^{(k)*}(a_{i})$ wherefrom
$\beta _{i} ~\Box  {}^{(k)}_{i}=0$, i.e., the constructed 
$\psi ^{(k+1)}$ is also equivalent
to $\varphi $.

To completely prove the Theorem it remains to make
use of the fact that $I^{k}_{{\cal U}}=0$ for a sufficiently great $k$.
\end{proof}

\begin{Corollary}If families $\varphi _{i}: {\cal U}
\tto \tilde \Cee ^{n|n}, \; i=1, 2$, are such
that $\varphi ^{*}_{1}(\tau_{i})=\varphi ^{*}_{2}(\tau_{i})$ for $i=1, \dots
, 2n$ and at least one of
the conditions
$$
\begin{matrix}
&(i) \cr
&(ii) \cr
&(iii) \end{matrix} \qquad \qquad
\begin{array}{ll}
&\varphi ^{*}_{1}(\alpha _{i})\neq 0 \cr
&\varphi ^{*}_{2}(\alpha _{i})\neq 0 \cr
&\varphi ^{*}_{1}(s_{i})\equiv \varphi ^{*}_{2}(s_{i}) (\mod  
I_{{\cal U}}) 
\end{array}
\qquad \qquad
\begin{matrix}
&for\; i=1, \dots, n \cr
&for\; i=1, \dots, n \cr
&for\; i=1, \dots, n \end{matrix}
$$
 is satisfied then the values of all the invariant
functions on $\varphi _{1}$ and $\varphi _{2}$ coincide.
\end{Corollary}

\begin{proof} Suppose either (i) or (iii) are satisfied.
Then by Theorem ?? applied to $\varphi =\varphi _{1}$ and $g_{i}=
\varphi ^{*}_{2}(s_{i})$ there exists
a $\psi $ equivalent to $\varphi _{1}$ such that 
$\varphi ^{*}(s_{i})=\varphi ^{*}_{2}(s_{i})$. If
$H=H(s_{1}, \dots
, s_{n}, \tau_{1}, \dots
, \tau_{n})$ is an invariant function then
$$
\begin{array}{ll}
\varphi ^{*}_{1}(H)&=\psi ^{*}(H)=H(\psi ^{*}(s_{1}), \dots
, \psi ^{*}(s_{n}), \varphi ^{*}_{1}(\tau_{1}), \dots
, \varphi ^{*}_{1}(\tau_{n}))= \cr
&=H(\varphi ^{*}_{2}(s_{1}), \dots
, \varphi ^{*}_{2}(s_{n}), \varphi ^{*}_{2}(\tau_{1}), \dots
, \varphi ^{*}_{2}(\tau_{n}))=\varphi ^{*}_{2}(H). 
\end{array}
$$
 If (ii) holds, then we have to interchange $\varphi _{1}$ and 
$\varphi _{2}$.
\end{proof}

\begin{Remark} The condition $g_{j}=\varphi ^{*}(s_{j}) 
(\mod  I_{{\cal U}})$ in the
theorem is essential: if some of $\varphi ^{*}(\alpha _{i})$ are equal to zero
then a set $b_{1}, \dots
, b_{n}$ for which $g_{j}=s_{j}(b_{1}, \dots
, b_{n}) (\mod  I_{{\cal U}})$
may not satisfy the condition $\prod \limits_{i<j}(b_{i}-b_{j})\neq 0 
(\mod  I_{{\cal U}})$ and
therefore further steps in the construction of $\psi ^{(s)}$
may prove impossible to perform. For instance, for $n=2, 
{\cal U}=\Cee ^{0|2}$ with coordinates $\xi _{1}, \xi _{2}$ 
and $\varphi ^{*}(\alpha _{i})=0, \varphi ^{*}(a_{1})=1, 
\varphi ^{*}(a_{2})=-1$ any pair of even functions $g_{1}, g_{2}$ satisfies
equations (3.5) but for $g_{1}=-2, g_{2}=1+\xi _{1}\xi _{2}$ there are no
functions $b_{1}, b_{2}$ on ${\cal U}$ such that $b_{1}+b_{2}=2, 
b_{1}b_{2}=1+\xi _{1}\xi _{2}$.
\end{Remark}

However, in the example considered, for calculation of
the values of invariants the functions $g_{1}$
and $g_{2}$ can be used instead of $s_{1}$ and $s_{2}$: they will vanish, 
anyway.

\section* {\protect \S 4. Description of invariant functions on $\Q(n)$ and 
$\Odd ~(n)$}

\ssec{4.1} Invariants of $\Q(n)$ and of $\Odd (n)$ have a
uniform description and, therefore, we will introduce the
following notations: $(M (n), \tilde M(n), G)$ stands for either of the
sets $(\Q(n)$, $\tilde \Q(n)$, $\GQ(n))$ or $(\Odd (n), 
\widetilde{\Odd} (n), \GL
(n|n))$.

On $M(n)$, there exists a set of $G$-invariant
polynomials $\tau_{k}$, where $k\in \Nee$:
$$
\tau_{k}(A)=\left\{\begin{matrix}k^{-1}\qtr  
A^{k}\; \text{ for }\; A\in \Q(n), \cr 
(2k-1)^{-1}\str A^{2k-1}\; \text{ for }\; A\in \Odd~ (n)\end{matrix},\right . 
$$
where $A$ is a family of matrices. 

Fix an embedding $j: \tilde \Cee ^{n|n}\tto 
\tilde M(n)$ having identified $\tilde \Cee ^{n|n}$ with the supermanifold of
(nonhomogeneous) diagonal matrices in $\tilde \Q(n)$ or block
matrices of the form $\begin{pmatrix} \alpha & A \cr 1_{n} & 0 \end{pmatrix}$ 
in
$\widetilde{\Odd} (n)$, where $A=\diag (a_{1}, \dots
, a_{n}), \alpha =(\alpha _{1}, \dots , \alpha _{n})$. The definition implies
that $j^{*}(\tau_{i})=\tau_{i}$, where in the left-hand side there
stands $\tau_{i}\in  F(\tilde M(n))$ and in
the right-hand side there
stands $\tau_{i}\in F(\tilde \Cee ^{n|n})$.

The embedding $j$ is compatible with the
$S_{n}\vdash \Cee ^{0|n}$-action on $\tilde \Cee ^{n|n}$  and the $G$-action on
$M(n)$ in the following sense.

\begin{Lemma} Two families of morphisms $\varphi _{i}: {\cal U}
\tto 
\tilde \Cee ^{n|n}, \; \quad i=1, 2$ pass into each other under the action of
$S_{n}\vdash 
\Cee ^{0|n}$ if and only if the families $j\circ \varphi _{1}$ and $j\circ 
\varphi _{2}$ pass into each other under the action of $G$.
\end{Lemma}

\begin{proof} Corollary 1.4 and arguments in 2.1 and 2.2
imply that the collection of eigenvalues of the family of
matrices from $\tilde M(n)$ is defined uniquely up to the
$S_{n}\vdash \Cee ^{0|n}$-action and it only remains 
to verify that the
equivalence of $\varphi _{1}$ and $\varphi _{2}$ implies the 
equivalence of $j\circ \varphi _{1}$
and $j\circ \varphi _{2}$.

Permutations are realized in a usual way and it
suffices to consider the action of $\Cee ^{0|n}$ on blocks: 
$$
(1+\varepsilon)(a+\alpha )(1+\varepsilon)^{-1}=a+\alpha +2\varepsilon \alpha 
\; \text{for $\Q(n)$}
$$
 and
$$
\begin{pmatrix} a+\varepsilon \alpha & \varepsilon a \cr 
\varepsilon & a \end{pmatrix}
\begin{pmatrix} 
\alpha & a \cr 1 & 0 \end{pmatrix}
\begin{pmatrix} a+\varepsilon \alpha & \varepsilon a \
\varepsilon & a\end{pmatrix} ^{-1}=
\begin{pmatrix} \alpha & a+2\varepsilon \alpha \\
1 & 0 \end{pmatrix}\text{  for  }\Odd (n).
$$
\end{proof} 

\begin{Remark} It would have been more natural to embed
$S_{n}\vdash \Cee ^{0|n}$ into $G$ in order to ensure that
this embedding commutes with $j$. It is clear, however, that such an
embedding exists for $\GQ(n)$ and does not exist for $\GL (n|n)$.
\end{Remark}

\begin{Theorem} For any point $m\in \tilde M(n)_{rd}$ there exist 
supermanifold morphisms $g_{m}: U_{m}\tto G$ and $\pi _{m}: 
U_{m}\tto \tilde \Cee ^{n|n}$ defined in a neighborhood of 
$m$ such that $j\circ \pi _{m}=\ad g_{m}: U\tto M(n)$ and 
the set of functions $\pi ^{*}_{m}(\alpha _{1}), \dots , \pi 
^{*}_{m}(\alpha _{n})$ can be complemented to a coordinate system on 
$U_{m}$ and the ideal generated by $\pi ^{*}(\alpha _{1}), \dots , \pi 
^{*}(\alpha _{n})$ coincides with the ideal generated by 
$\tau_{1}|_{U_{m}}, \dots , \tau_{n}|_{U_{m}}$.
\end{Theorem}

{\bf Proof}~ for $M=\Odd~ (n)$ reduces to the fact that having fixed a 
basis $e_{1}, \dots , e_{n}$, $\varepsilon _{1}, \dots , \varepsilon 
_{n}$ in a free $n|n$-dimensional module $L$ over $F(U)$ we can 
consider a neighborhood $U_{m}$ of point $m$ as a linear operator on 
$L$ that we will denote by $A$.

As we will show in what follows, in a neighborhood of $m$ there are 
defined projections $P_{i}: L\tto L, \; i=1, \dots , n$, to 
$A$-invariant submodules.  Let us select an even vector $y$ and an odd 
vector $\eta $ in $L$ such that the set $P_{1}y, \dots , P_{n}y, 
P_{1}\eta , \dots , P_{n}\eta $ is a basis in $L$.  If $m\in \tilde 
\Cee ^{n|n}_{rd}$ then for $y$ and $\eta $ we can take $e_{1}+\dots 
+e_{n}$ and $\varepsilon _{1}+\dots +\varepsilon _{n}$ and in the 
general case set $y=g^{-1}(e_{1}+\dots +e_{n}), \eta 
=g^{-1}(\varepsilon _{1}+\dots +\varepsilon _{n})$, where $g\in 
G_{rd}$ is such that $gmg^{-1}\in \tilde \Cee ^{n|n}_{rd}$.

The pairs $P_{i}y, P_{i}\eta $ constitute bases of $A$-invariant
submodules and, therefore, for every $i$ we have a morphism
$U_{m}\tto \Cee ^{1|1}$ which gives rise to a morphism 
$\pi _{m}:  U_{m}\tto \tilde \Cee ^{n|n}$ and the transition matrix
from the basis $e_{1}, \dots , e_{n}, \varepsilon _{1}, \dots , \varepsilon
_{n}$ to the basis $P_{1}y, \dots , P_{n}y, P_{1}\eta , \dots , P_{n}\eta $ is
the desired $g_{m}: U_{m}\tto G$.

For $M=\Q(n)$ the proof differs only in that $L$ is a free right 
module of $\rk~ \; n$ and we select one vector $y$ such that $P_{1}y, 
\dots , P_{n}y$ is a basis in $L$.

It only remains to prove the existence of projections $P_{i}$.  Let 
$\mu $ be an eigenvalue of a complex matrix $m\in \tilde M(n)_{rd}, 
V\subset \Cee $ an open disk whose interior contains $\mu $ and does 
not contain eigenvalues of $m$ distinct from $\mu $.  Then any matrix 
$m'$ from a neighborhood $U_{rd}\ni m$ has only one eigenvalue $\mu '$ 
in $V$ (we diminish $U$ if necessary without much ado).

Similarly to \cite{RS}, set
$$
P(A)=\left\{
\begin{matrix}
&\frac{1}{2\pi i}\int _{\partial V}(\lambda E_{n}-A)^{-1}d\lambda 
\; \text{for $M(n)=\Q(n)$} \cr
&\frac{1}{2\pi i}\int _{\partial V}(\lambda E_{n|n}-A^{2})^{-1}d
\lambda \; \text{for $M(n)=\Odd (n)$}\end{matrix}.
\right.
$$
Clearly, $P$ is an even operator commuting with $A$.  Let us establish 
that $P$ is indeed a projection to a 1-dimensional submodule if 
$M(n)=\Q(n)$ or $1|1$-dimensional if $M(n)=\Odd (n)$.

On $M(n)$, there is a standard global coordinate system which 
determines the factorization $U\equiv U_{rd}\times \Cee ^{0|k}$, where 
$k$ is the number of odd coordinates.  Therefore, it is possible to 
assume that the matrices of operators $P$ and A are $U_{rd}\times \Cee 
^{0|k}$-families of matrices.  For any point $m'\in U_{rd}$ the 
corresponding $\Cee ^{0|k}$-families of matrices $A(m')$ and 
$P(A(m'))=P((A)m')$ are matrices with elements from the finite 
dimensional Grassmann algebra $\Lambda =F(\Cee ^{0|k})$ and, 
therefore, there exists an even invertible matrix $g$ with elements 
from $\Lambda $ such that $g \cdot A(m') \cdot g^{-1}$ is of the 
standard format; hence, $g\cdot P((A)(m'))\cdot g^{-1}=P(g\cdot A\cdot 
g^{-1}(m'))$ and coincides with the projection onto the submodule 
corresponding to eigenvalue $\mu '$.

Actually, $P(A)$ is \lq\lq composed" from the projections
$P(m')$ corresponding to eigenvalues $\mu '$ and existing for
every $m'\in U_{rd}$ and the explicit formula for $P(A)$
establishes a holomorphic dependence of $P(m')$.

Since $\tau_{k}$ are invariant functions, then
$$
\tau_{k}|_{U_{m}}=\pi ^{*}_{m}(\tau_{k})=\sum \pi ^{*}_{m}
(a^{k-1}_{i})\cdot \pi ^{*}_{m}(\alpha _{i})
$$
and, therefore, the passage from the collection $\tau_{1}|_{U_{m}}, 
\dots , \tau_{n}|_{U_{m}}$ to $\pi ^{*}_{m}(\alpha _{1}), \dots , \pi 
^{*}_{m}(\alpha _{n})$ is performed by an invertible linear 
transformation.  By a $G$-action an arbitrary point of $M(n)_{rd}$ can 
be transformed into a point of $j(\tilde \Cee ^{n|n})_{rd}$ and, 
clearly, in a neighborhood of $j(\tilde \Cee ^{n|n}_{rd})$ the 
functions $\tau_{1}, \dots , \tau_{n}$ can be complemented to a local 
coordinate system.

In what follows the pairs of morphisms $\pi : U\tto \tilde 
\Cee ^{n|n}, g: U\tto G$ possessing properties established 
in the theorem will be called {\it projections}.

\ssec{4.2} We intend to establish an isomorphism between $G$-invariant 
functions on $\tilde M(n)$ and $S_{n}\vdash \Cee ^{0|n}$-invariant 
functions on $\tilde \Cee ^{n|n}$ and, therefore, with balanced 
functions on $\check \Cee ^{n|n}$.  For this it is necessary to lift 
the functions $s_{i}$ to $\tilde M(n)$.

\begin{Theorem} There exist even rational functions $s_{1}, \dots , 
s_{n}$ on $M(n)$ without singular points on $\tilde M(n)$ and 
satisfying the system of equations
$$
\tau_{k+n+1}=\sum^{}_{1\le i\le h}\tau_{k+i}s_{n-i+1}, \; \; k=0, 
\dots , n-1 \eqno{(4.1)}
$$
For $n>1$ there are no polynomial solutions of system (4.1).
\end{Theorem}

{\bf Proof of Theorem}.  Let the $x_{i}, \xi _{j}$ be standard 
coordinates on $M(n)$.  Let us express $\tau_{1}, \dots , \tau_{2n}$ 
and the functions $s_{1}, \dots , s_{n}$ to be described in the form
$$
\tau_{i}=\sum^{}_{\alpha }c^{\alpha }_{i}(x)\xi ^{\alpha }, s_{i}=
\sum^{}_{\beta }d^{\beta }_{i}(x)\xi ^{\beta }, 
$$
where $\alpha $ and $\beta $ run over sets of $0$'s and $1$'s of 
length $k$, $c^{\alpha }_{i}$ are known polynomials and $d^{\beta 
}_{i}$ unknown functions.  Then equating coefficients of $\xi ^{\alpha 
}$ in the left and right-hand sides of (4.1) we get an equivalent to 
(4.1) system of linear non-homogeneous equations in functions 
$d^{\beta }_{i}$, where the coefficients and constant terms are 
polynomials on $M_{rd}$.  Let us call this system {\it the main} one 
but will not write it.

In order to avoid confusion, let us denote for the
time being the functions $\tau_{i}$ and $s_{i}$ on $\tilde 
\Cee ^{n|n}$ by $\tau'_{i}$ and
$s'_{i}$.

Let $m\in M(n)_{rd}$, $\pi _{m}: U_{m}\tto \tilde \Cee 
^{n|n}$; let $ g_{m}: U_{m}\tto G$ be the projection that 
exists by Lemma.  Then $\pi ^{*}(\tau'_{i})=\pi ^{*} \circ 
j^{*}(\tau_{i})=\tau_{i}$ since $j\circ \pi _{m}=\ad g_{m}$.  The 
functions $\pi ^{*}(s'_{1}), \dots , \pi ^{*}(s'_{n})$ form a solution 
of the main system on $U$ since (4.1) on $\tilde \Cee ^{n|n}$ is 
identically satisfied.

Therefore, the main system is compatible in a neighborhood of any 
point of $M(n)_{rd}$.  Since its coefficients are polynomials then for 
any point $m\in M(n)_{rd}$ there exists a solution that can be 
extended to a Zariski open neighborhood, a solution that consists of 
rational functions on $M(n)_{rd}$ and has no singularities at $m$.  
The sheaf ${\cal P}$ of solutions of the system of homogeneous 
equations corresponding to the main system is coherent and $\tilde 
M(n)_{rd}$ is an affine algebraic variety (singled out in $M(n)_{rd}$ 
by the condition $f(m)=0$, where $f$ is the discriminant of the 
characteristic polynomial of $m)$ and therefore $H^{1}(\tilde 
M(n)_{rd}, {\cal P})=0$ by Serre's theorem.  This means that there 
exists a global solution of the main system --- the set of rational 
functions $d^{\beta }_{i}$ without singularities on $\tilde 
M(n)_{rd}$.

Setting $s_{i}=\sum d^{\beta }_{i}\xi ^{\beta }$ we get the required 
solution of system (4.1).

{\bf Remarks}. 1) At the moment we cannot explicitly produce the
functions $s_{i}$.

2) The set of functions $s_{i}$ is by no means unique but
in what follows we will fix one such set.

If $h_{1}, \dots , h_{n}$ is a solution of system (4.1) and $n>1$, 
then $h_{1}$ satisfies the equation
$$
h_{1}\tau_{1}\cdot \dots
\cdot \tau_{n}=\tau_{1}\cdot \dots
\cdot \tau_{n-1}\tau_{n+1}
\eqno{(4.2)}
$$

Let us consider $\Q(n)$ and $\Odd (n)$ separately.  On $Q (n)$ the 
even and odd coordinates fill in two square matrices $B$ and $\beta $, 
respectively, and $\tau_{k}= k^{-1}\qtr (B+\beta )^{k}$ is an 
homogeneous polynomial of degree $k$ in $B$ and $\beta $ such that 
$\beta $ is only encountered in odd degrees and the highest with 
respect to $B$ and (simultaneously) lowest with respect to $\beta $ 
term is equal to $\tr B^{k-1}\beta $; the second highest in $B$ term 
is of degree $k-3$.

The functions $\tau_{1}, \dots , \tau_{n}$ can be included into a 
local coordinate system on $\Q(n)$ and, therefore, the degree of the 
product $\tau_{1}\dots \tau_{n}$ with respect to $\beta $ is equal to 
$n$ and the highest in $B$ term in $\tau_{1}\dots \tau_{n}$ is 
$\prod\limits _{1\le i\le n}\tr B^{i-1}\beta $.  After the change 
$g_{1}=\Delta _{1}+\tr B$ the equation (4.2) turns into
$$
\Delta _{1}\tau_{1}\dots
\tau_{n}=\tau_{1}\dots
\tau_{n-1}(\tau_{n+1}-\tr B\cdot \tau_{n})
\eqno{(4.3)}
$$
The degree of $\tau_{1}\cdot \dots \cdot \tau_{n-1}(\tau_{n+1}-\tr 
B\cdot \tau_{n})$ with respect to $B$ does not exceed the degree of 
$\tau_{1}\cdot \dots \cdot \tau_{n-1}(\tau_{n+1}-\tr B^{n}\cdot \beta 
)$ which is less than the degree of $\tau_{1}\cdot \dots \cdot 
\tau_{n}$.  The point is that the highest with respect to $B$ term in 
$\tr B^{n}\cdot \beta -\tr B\cdot \tau_{n}$ is equal to $\sum\limits 
^{}_{1<i<n-1}k_{i}(B) \tr (B^{i-1}\beta )$, where $k_{i}(B)$ are 
polynomials and the summand $k_{i}(B)\tr (B^{i-1} \beta )$ is killed 
being multiplied by the highest term of $\tau_{i}$.  Therefore, the 
main system has no polynomial solutions for $\Q(n)$.

For $\widetilde{\Odd} (n)$ the coordinates fill in the matrix 
$\begin{pmatrix} \alpha & B \cr C & \delta \end{pmatrix}$, where 
$\alpha $ and $\delta $ consist of odd coordinates and where $B$ and 
$C$ consist of even coordinates.

Set ${\cal V} =\begin{pmatrix} \alpha & 0 \cr 0 & \delta 
\end{pmatrix}, {\cal U}=\begin{pmatrix} 0 & B \cr C & 0 \end{pmatrix}$ 
.  Then $\tr {\cal U}^{2k-2}{\cal V} $ is the lowest term of 
$\tau_{k}=\frac{1}{2k-1}\str ({\cal V} + {\cal U})^{2k-1}$ with 
respect to ${\cal V} $ and simultaneously the highest term with 
respect to ${\cal U}$.  The change $g_{1}=\tr BC+\Delta $ turns (4.2) 
into
$$
\Delta \tau_{1}\dots
\tau_{n}=\tau_{1}\dots
\tau_{n-1}(\tau_{n+1}-\tr BC\cdot \tau_{n}).
\eqno{(4.4)}
$$
In exactly the same way as this was done for $\Q(n)$ it is easy to 
show that for $n>1$ the degree of the right-hand side of (4.4) with 
respect to ${\cal U}$ is smaller than the degree of $\tau_{1}\dots 
\tau_{n}$ and there are no polynomial solutions of the main system.

\ssec{4.3} Clearly, an infinite dimensional supergroup of morphisms 
from $U$ to $G$ acts on the set of solutions (4.1) defined on 
$U\subseteq M(n)$.  This completely describes the nonuniqueness of the 
set $s_{1}, \dots , s_{n}$.

\begin{Theorem} If $f_{1}, \dots
, f_{n}$ are even functions on an open
subsupermanifold $W\subset \tilde M(n)$ satisfying (4.1), then

$1)$ For any point $m\in U_{rd}$ there exists a morphism $h_{m}: 
V\tto G$ defined in a neighborhood of $m$ that sends the 
set $s_{1}, \dots , s_{n}$ to $f_{1}, \dots , f_{n}$;

$2)$ The functions $\cpr  f_{i}$ on $U_{rd}$ are determined
from the relation
$$
\lambda ^{n}+\sum^{}_{0\le i\le n-1}\lambda ^{i}\cpr  
f_{n-i}(m)=\det (\lambda E-A(m)), 
$$
 where
$$
A(m)=\left\{\begin{matrix} m \; \text{ for}\; m\in \tilde \Q(n)_{rd} 
\cr BC \; \text{ for}\; m=\begin{pmatrix} 0 & B \cr C & 0\end{pmatrix} 
\in \widetilde{\Odd} (n)_{rd} \end{matrix}\right..
$$
\end{Theorem}

\begin{proof} Let $g_{m}: U\tto G$ and $\pi _{m}: 
U\tto \tilde \Cee ^{n|n}$ be the projection (see 4.1).  
Then having considered $\pi _{m}$ as a $U$-family of points $\tilde 
\Cee ^{n|n}$ let us apply Theorem 3.4 to the collection of functions 
$f_{1}, \dots , f_{n}$: they satisfy (3.5) and $\pi ^{*}_{m}(\alpha 
_{i})\neq 0$ and, therefore, $f_{i}\equiv \pi ^{*}(s_{i})(\mod I)$, 
where $I$ is the ideal generated by odd coordinates on $U$ and there 
exists a family of morphisms $\pi '_{m}: U_{m}\tto \tilde 
\Cee ^{n|n}$ equivalent to $\pi _{m}$ such that $f_{i}=(\pi 
'_{m})^{*}(s_{i})$.

Since functions $\pi ^{*}_{m}(\alpha _{i})$ can be included into a 
local coordinate system and equivalence of $\pi _{m}$ and $\pi '_{m}$ 
means that
$$
\pi ^{*}_{m}(\alpha _{i})=\pi '{*}(\alpha _{i}) \; {\rm and}\; 
\pi ^{*}_{m}(\alpha _{i})[\pi ^{*}_{m}(a_{i})-\pi '{*}_{m}(a_{i})]=0, 
$$
then there exist odd functions $k_{1}, \dots , k_{n}$ such that $\pi 
'{*}_{m}(a_{i})=\pi ^{*}_{m}(a_{i})-\pi ^{*}_{m} (\alpha _{i})k_{i}$.  
In other words there exists a $U$-family of points of $S_{n}\vdash 
\Cee ^{0|n}$ that sends $\pi '_{m}$ to $\pi _{m}$.  Having identified 
$\tilde \Cee ^{n|n}$ with $j(\tilde \Cee ^{n|n})$ let us lift this 
family to a $U$-family of points of $G$ that sends $j\circ \pi '_{m}$ 
to $j\circ \pi _{m}$ and, therefore, sends the set $\pi 
^{*}_{m}(s_{i})$ to the set $f_{i}$.  Thus, the sets $\{f_{i}\}$ and 
$\{s_{i}\}$ can be locally obtained from the set $\{\pi 
^{*}_{m}(s_{i})\}$, hence from each other.

We have proved above that $\cpr f_{i}=\cpr s_{i}$ and therefore the 
functions $\cpr s_{i}$ are $G_{rd}$-invariant.  The formulas we are 
proving for $\cpr s_{i}$ are also $G_{rd}$-invariant; they are 
satisfied on $j(\tilde \Cee ^{n|n})_{rd}$ and therefore they are true 
on $\tilde M(n)$.  \end{proof}

\ssbegin{4.4}{Theorem} The algebra of invariant functions on $\tilde 
M(n)$ is isomorphic to the algebra of $S_{n}\vdash \Cee 
^{0|n}$-invariant functions on $\tilde \Cee ^{n|n}$ and, therefore, to 
the algebra of balanced functions and the isomorphism is performed via 
$j^{*}: F(\tilde M(n))\tto F(\tilde \Cee ^{n|n})$ under which the 
polynomials on $\tilde M(n)$ are identified with polynomials on 
$\tilde \Cee ^{n|n}$ and rational functions with rational functions.
\end{Theorem}

\begin{proof} The properties of $j^{*}$ immediately imply that if $f$ 
is an invariant function on $\tilde M(n)$ then $j^{*}(f)$ is an 
invariant function on $\tilde \Cee ^{n|n}$ and since locally $f$ 
coincides with $\pi ^{*}_{m}\circ j^{*}(f)$ then $j^{*}(f)=0$ implies 
$f=0$.  Now let $f'$ be an invariant function on $\tilde \Cee ^{n|n}$ 
and $\check f$ the corresponding balanced function on $\check \Cee 
^{n|n}$.

In a neighborhood of any point $m\in M(n)_{rd}$ we can apply Theorem 
3.4 to the $U$-family $\pi _{m}: U\tto\tilde \Cee ^{n|n}$ 
and deduce that the functions $\tau_{1}, \dots , \tau_{n}$ are 
$G$-invariant, the $G$-action sends $s_{1}, \dots , s_{n}$ to another 
solution of (4.1) and, therefore, the function $\check f(s_{1}, \dots 
, s_{n}, \tau_{1}, \dots , \tau_{n})$ does not vary, i.e., $\check 
f(s_{1}, \dots , s_{n}, \tau_{1}, \dots , \tau_{n})$ is an invariant 
function.

If $f'$ is an invariant polynomial then by Theorem 3.3 $f'=P(\tau_{1}, 
\dots , \tau_{k})$ is an invariant polynomial in $\tau_{1}, \dots , 
\tau_{k}$ on $\tilde \Cee ^{n|n}$; hence, $P(\tau_{1}, \dots , 
\tau_{k})$ is an invariant polynomial on $\tilde M(n)$.  If $f'$ is a 
rational function then $f$ is also a rational function since $\check 
f$ and $s_{1}, \dots , s_{n}$ are rational functions.
\end{proof}

\begin{Corollary} Any invariant polynomial $P$ on $M(n)$ can be
uniquely expressed in the form
$$
P=\sum^{}_{0\le k\le n}\; \; \sum^{}_{i_{1}<\dots
<i_{k}}c^{i_{1}\dots
i_{k}}\tau_{i_{1}}\dots
\tau_{i_{k}}
$$
where only a finite number of coefficients $c^{i\dots i}\in \Cee $ is 
nonzero.  All the relations between the $\tau_{1}, \dots , \tau_{n}, 
\dots $ are corollaries of supercommutativity and relations 
$\tau_{i_{1}}\dots \tau_{i_{n+1}}=0$.
\end{Corollary}

\ssec{4.5.  The case of $\tilde M(n)$} In exactly the same way as for 
$\tilde \Cee ^{n|n}$, to compute the values of any function it 
suffices to know $\tau_{1}, \dots , \tau_{2n}$.

\begin{Theorem} Let $\varphi : {\cal U}\tto \tilde M(n)$ be 
a family of points of $\tilde M(n)$ and $h_{1}, \dots , h_{n}$ even 
functions on ${\cal U}$ satisfying the equations
$$
\varphi ^{*}(\tau_{n+1+k})=\sum^{}_{1\le i\le k}
\varphi ^{*}(\tau_{n+k+1-i})h_{i}, \; \text{where } k=0, \dots
, n-1\text{ and }h_{i}\equiv \varphi ^{*}(s_{i}) 
\pmod I_{{\cal U}}).
\eqno{(4.5)}
$$
Then for any invariant function $f$ on $\tilde M(n)$ we have $\varphi 
^{*}(f)=\check f(h_{1}, \dots , h_{n}, \varphi ^{*}(\tau_{1}), \dots , 
\varphi ^{*}(\tau_{n}))$, where $\check f$ is the balanced function 
corresponding to $f$.
\end{Theorem}

\begin{proof} Let us apply Theorem 3.4 to the family
$\pi _{\varphi }\circ \varphi : V\tto \tilde 
\Cee ^{n|n}$ defined in a neighborhood of a 
point $u\in {\cal U}_{rd}$.
\end{proof}

\begin{Corollary} If the first $2n$ invariant polynomials of
the families of morphisms $\varphi _{i}: {\cal U}
\tto \tilde M(n), \; \; i=1, 2$, coincide and
$\varphi ^{*}_{1}(s_{i})\equiv 
\varphi ^{*}_{2}(s_{i}) (\mod  I_{{\cal U}})$ then 
the remaining invariant
functions also coincide.
\end{Corollary}

\begin{Remark} 1) Recall that $\cpr  (s_{i})$ are invariant
polynomials on $M(n)_{rd}$ which allows to solve system (4.5).

2) It seems strange that it is possible to determine the value of any 
invariant function $f$ from $\tau_{1}, \dots , \tau_{2n}$ whereas $f$ 
cannot as a rule be expressed in the form of a function in $2n$ odd 
variables $\tau_{1}, \dots , \tau_{2n}$.  The point is that the 
variables $\tau_{1}, \dots , \tau_{2n}$ are not independent: the 
product of any $n+1$ of them is equal to zero.
\end{Remark}

\ssec{4.6} The collection of functions $s_{1}, \dots , s_{n}$ not only 
gives a collection of invariants with values in $\Cee $ but is a set 
of $G$-invariants of \lq\lq the second turn" in the following precise 
sense.

Denote by $L(n)$ the closed subsupermanifold in $\tilde M(n)$
singled out by equations 
$$
\tau_{1}=0, \dots
, \tau_{n}=0.
$$ 
It is $G$-invariant
together with $\tau_{i}$ and $L(n)_{rd}=\tilde M(n)_{rd}$.

\begin{Theorem} {\em 1)} If $f$ is an invariant function on $\tilde M(n)$
then $f|_{L(n)}$ is a constant.

$2)$ The functions $l_{1}=s_{1}|_{L(n)}, \dots , l_{n}=s_{n}|L(n)$ do 
not depend on the choice of the collection $s_{1}, \dots , s_{n}$ --- 
a solution of (4.1) --- and are generators of the algebra of invariant 
functions on $L(n)$.
\end{Theorem}

{\bf Proof}. 1) Lemma (3.3) implies that $\check f\equiv {\rm const} 
(\mod \tau_{1}, \dots , \tau_{n})$.

2) If $\pi _{i}: U\tto \tilde \Cee ^{n|n}, \; 
\; i=1, 2$, are two projections then on $L(n)$ 
$$
\pi ^{*}_{1}(\alpha _{i})=\pi ^{*}_{2}(\alpha _{i})=0\; {\rm and 
}\; \pi ^{*}_{1}(a_{i})=\pi ^{*}_{2}(a_{\delta (i)}), \; {\rm where 
}\; \delta \in S_{n}, 
$$ 
and, therefore, all solutions $s'_{1}, \dots , s'_{n}$ of (4.1) give 
the same set of functions $l_{1}, \dots , l_{n}$ which are basis 
symmetric functions in $\pi ^{*}(a_{1}), \dots , \pi ^{*}(a_{n})$ 
defined uniquely up to the $S_{n}$-action.  \qed

\section* {\protect \S 5. Examples}

\ssec{5.1} On the image of the embedding$j(\tilde \Cee ^{n|n})
\hookrightarrow \Q(n)$ 
we have $\qet  A=\sum \alpha _{i}a^{-1}_{i}$. In
particular, for $n=2$ we have
$$
\qet  A=\frac{\tau_{1}(a_{1}+a_{2})-\tau_{2} }{ a_{1}a_{2}}.
$$
If the family of matrices $\psi : {\cal U}\tto \Q(n)$ is 
such that the functions $\tau_{l}=l^{-1}\qtr A^{l}$ are defined for 
$l$ close to 0, then $\qet~A =\lim\limits_{l\to 0}\tau_{l}$.

On $\widetilde{\Odd} (n)$, the expression $\sum \alpha _{i}a^{-1}_{i}$ 
also defines an invariant function, $\tau_{0}=-\str A^{-1}$, that does 
not possess, unlike qet, any special properties; in particular, 
$\tau_{0}(kA)=k^{-1}A$.  However, both $\qet (\lambda -A)$ on $\Q(n)$ 
and $-\str (\lambda -A)^{-1}$ on $\Odd (n)$ are generating functions 
for all invariants of $A$:
$$
\begin{matrix}
\qet (\lambda -A)&=\sum^{}_{i}\frac{d_{i}}{a_{i}-\lambda }=
-\sum^{\infty }_{j=0}\lambda ^{-j}\tau_{j}(A); \cr
-\str  (\lambda -A)^{-1}&=-\lambda ^{-1}\str 
(1-\lambda ^{-1}A)^{-1}= \cr
&=-\sum^{\infty }_{j=0}\lambda ^{-(j+1)}\str A^{j}=
-\sum^{\infty }_{j=0}\lambda ^{-(2j+2)}\str A^{2j+1}. \end{matrix}
$$

\ssec{5.2} On $\Q(2)$, with coordinates that fill in two
square matrices: an even one, $B=(b_{ij})$, and an odd one, 
$\beta =(\beta _{ij})$, one of the
rational solutions of (4.1) is given by the formulas
$$
\begin{matrix}
s_{1}(B, \beta )&=b_{11}+b_{22}+2\frac{(\beta _{22}-\beta _{11}) 
[\beta _{12}b_{21}-\beta _{21}b_{12}]+(b_{11}-b_{22})\beta _{12}\beta 
_{21} }{ (b_{11}-b_{22})^{2}+4b_{12}b_{21}} \cr 
s_{2}(B, \beta 
)&=\frac{1 }{ 2}s_{1}(B+\beta ^{2}, B\beta +\beta B).
\end{matrix}
$$

\ssec{5.3} The results of \S 4 provide us with a complete system of 
invariants for a linear superbundle of rank $n|n$ in each fiber of 
which there is fixed an odd invertible operator $A$ --- the set of 
polynomials
$$
\str A, \; \dots
, \; \frac{1 }{ 4n-1}\str A^{4n-1}.
$$
Two particular cases are of special interest: (a) one, connected with 
an almost complex structure \cite{Po}, \cite{Va} and (b) another, 
connected with a pair (a symplectic structure, a periplectic 
structure), cf.  \cite{Kh}, \cite{V}, \cite{KN}, \cite{NK}.

\ssec{5.3.1} In the first case $A^{2}=-1_{n|n}$, i.e., $A$ is \lq\lq 
far" from the general position.  As shown in \cite{Po}, \cite{Va} such 
operators have no invariants: by a change of basis $A$ can be reduced 
to the form $\begin{pmatrix} 0 & -1_{n} \cr 1_{n} & 0 \end{pmatrix}$.

Let us show that a more general statement is also true: the family of 
odd invertible matrices $A$ reduces to the form $\begin{pmatrix} 0 & 
B_{1} \cr B_{2} & 0 \end{pmatrix}$ if and only if $A^{2}$ reduces to 
the form $\begin{pmatrix} C_{1} & 0 \cr 0 & C_{2}\end{pmatrix}$.

Indeed, let $\begin{pmatrix} X & Y \cr Z & T \end{pmatrix}^{2}= 
\begin{pmatrix} C_{1} & 0 \cr 0 & C_{2}\end{pmatrix}$.  Then by acting 
on $A$ with the matrix $\begin{pmatrix} 1 & 0 \cr 0 & Z\end{pmatrix}$ 
(by conjugations) we may make $Z=1$; then $X+T=0$ and
$$
\begin{pmatrix} 1 & -X \cr 0 & 1 \end{pmatrix}\begin{pmatrix} X & Y 
\cr 1 & -X
\end{pmatrix}
\begin{pmatrix} 1 & X \cr 0 & 1 \end{pmatrix}=\begin{pmatrix} 0 & 
Y+X^{2} \cr 1 & 0 \end{pmatrix}.
$$
Notice that the results of \S 4 are
inapplicable here since $A$ has only two eigenvalues: $\pm 1$.

\ssec{5.3.2} If on a $2n|2n$-dimensional supermanifold there are given 
an even and an odd closed nondegenerate differential 2-forms $\omega 
_{\bar 0}$ and $\omega _{\bar 1}$ then $A=\omega ^{-1}_{\bar 0} \omega 
_{\bar 1}$ is an invertible odd linear operator in the tangent bundle.  
If the pair $(\omega _{\bar 0}, \omega _{\bar 1})$ is in the general 
position then the eigenvalues of $\cpr A$ are distinct (for $n=1$ this 
holds automatically) and we can make use of the results of \S 4.  The 
skew-symmetry of $\omega _{\bar 0}$ and $\omega _{\bar 1}$ as bilinear 
forms leads to the fact that $\tau_{2k-1}(A)=0$ and, therefore, all 
the invariants of $\omega _{\bar 0}$ and $\omega _{\bar 1}$ that can 
be obtained from $\omega ^{-1}_{\bar 0}\omega _{\bar 1}$ are 
$\tau_{2}, \dots , \tau_{4n-2}$.

\section* {\protect \S 6. Conclusion}

In this work we have obtained the complete set of
invariant functions on $\Q(n)$ and $\Odd (n)$ and there is a
constructible recipe for computing the values of any
invariant function from $\tau_{1}, \dots
, \tau_{2n}$. Apart from this
concrete information certain more abstract considerations
seem to be useful.

It is natural to interpret results of \S 4 as
follows: the quotient manifold $\tilde M/G$ does not exist in the
category of manifolds but exists in a broader category of
virtual supermanifolds \cite{L}, where
$$
\tilde M/G=\tilde \Cee ^{n|n}/S_{n}\vdash \Cee ^{0|n}
$$
and $G$-invariant functions on $\tilde M$ are \lq\lq functions" 
(whatever this might mean) on $\tilde M/G$.  In such terms the 
mysterious problem of computing invariant functions in $\tau_{1}, 
\dots , \tau_{2n}$ with the help of intermediary non-uniquely defined 
functions $s_{1}, \dots , s_{n}$ and balanced functions on $\check 
\Cee^{n|n}$ means, it seems, that all the functions on $\tilde M/G$ 
can be expressed as functions on the virtual supermanifold 
distinguished in $\Cee ^{0|2n}$ by equations $\tau_{i_{1}}\dots 
\tau_{i_{n+1}}=0$ for any $i_{1}, \dots , i_{n+1}$.

At the moment there is no theory of virtual supermanifolds.  Though in 
\cite{L}, $\# 31$, there are given examples of virtual supermanifolds 
which are not supermanifolds, the virtual supermanifolds were mainly 
introduced as a convenient means of work with \lq\lq genuine" 
supermanifolds.  The results obtained above can be considered as an 
experimental data contributing to the theory of virtual 
supermanifolds.

\end{document}